\documentclass[del-ppn]{ppn}
\usepackage{graphics}
\usepackage{amssymb,amsfonts,amsmath,latexsym}
\newcommand{\AMS}{\noindent{\small\sl AMS classification (2000). }}
\newcommand{\eqn}[1]{(\ref{#1})}
\newcommand{\be}{\begin{equation}}
\newcommand{\ee}{\end{equation}}
\newcommand{\finprf}{\unskip\null\hfill$\;\square$\vskip 0.3cm}
\newcommand{\N}{\mathbb{N}}
\newcommand{\R}{\mathbb{R}}
\newcommand{\C}{\mathbb{C}}

\newcommand{\X}{X}
\newcommand{\Y}{Y}
\newcommand{\W}{Z}
\newcommand{\aatop}[2]{\genfrac{}{}{0pt}{}{#1}{#2}}
\newcommand{\llc}{\lambda_1^{\mathcal L}}
\newcommand{\lambdaL}{\lambda^{\mathcal L}}
\newcommand{\PiL}{\Pi_{\mathcal L}}
\newcommand{\PiLC}{\Pi_{\mathcal L}^c}
\newcommand{\PiTot}{\Pi}
\newcommand{\PiTotC}{\Pi^c}
\newcommand{\BL}{B_{\mathcal L}}
\newcommand{\dnu}{d_+(\nu)}
\newcommand{\deltanu}{\delta}
\newcommand{\nnrm}[1]{|\kern-1pt|\kern-1pt|#1|\kern-1pt|\kern-1pt|}

\begin{document}
\baselineskip=6pt

\title{Relativistic hydrogenic atoms\\
in strong magnetic fields}
\author{Jean Dolbeault\inst{1} \and Maria J. Esteban\inst{1} \and Michael Loss\inst{2}}
\institute{Ceremade (UMR CNRS no. 7534), Universit\'e Paris Dauphine, Place de Lattre de Tassigny, 75775 Paris C\'edex~16, France \and School of Mathematics, Georgia Institute of Technology, Atlanta, GA 30332, USA}
\date{\today}
%

\maketitle 

\begin{abstract} In the Dirac operator framework we characterize and estimate the ground state energy of relativistic hydrogenic atoms in a constant magnetic field and describe the asymptotic regime corresponding to a large field strength using relativistic Landau levels. We also define and estimate a critical magnetic field beyond which stability is lost. \end{abstract}

\par\medskip
\AMS{\scriptsize {35Q40, 35Q75, 46N50, 81Q10; 34L40, 35P05, 47A05, 47N50, 81V45}

}

\section{Introduction}\label{Sec:Intro}

In this paper we characterize the ground state energy of hydrogenic atoms in magnetic fields. We deal with fields of large strength in the Dirac operator framework, far away from the perturbative regime.

To compute eigenvalues of Dirac operators, the usual min-max principle does not apply. More sophisticated versions of this principle have been established over the last few years, see \cite{GS,GLS,DES}. These techniques are powerful enough to provide accurate and efficient algorithms for calculating eigenvalues of Dirac operators \cite{DESV,DES2}. In this paper we demonstrate that they are also flexible enough to cover the case with a magnetic field and provide reasonable results for a highly non-perturbative problem, when paired with the right physical insight.

\medskip The Dirac operator for a hydrogenic atom in the presence of a constant magnetic field $B$ in the $x_3$-direction is given by
\be
H_B- \frac{\nu}{|x|} \quad\mbox{with}\quad H_B:= \alpha \cdot \left[\frac{1}{i}\,\nabla + \frac{1}{2}\,B(-x_2,x_1, 0)\right] +\beta \; , \label{one}
\ee
where $\nu = Z \alpha<1$, $Z$ {is} the nuclear charge number. The Sommerfeld fine-structure constant {is} $\alpha \approx 1/137.037$. The energy is measured in units of $mc^2$, {\sl i.e.,\/} the rest energy of the electron, the length in units of {$\hbar/mc$}, {\sl i.e.,\/} the Compton wavelength divided by $2\,\pi$, and the magnetic field strength $B$ is measured in units of {$\frac{m^2 c^2}{|q| \hbar} \approx 4.4\times 10^9$ Tesla}. Here $m$ is the mass of the electron, $c$ the speed of light, $q$ the charge of the electron (measured in Coulomb) and $\hbar$ is Planck's constant divided by $2 \pi$. It is worth recalling the the earth's magnetic field is of the order of $1$ Gauss and $1$ Tesla is $10^4$ Gauss.

In (\ref{one}), $\alpha_1$, $\alpha_2$, $\alpha_3$ and $\beta$ are $4 \times 4$ complex matrices, whose standard form (in $2\times 2$ blocks) is
$$
\beta=\left( \begin{matrix} {
\mathbb I} & 0 \\ 0 & -{\mathbb I} \\ \end{matrix} \right),\quad\alpha_k=\left( \begin{matrix} 0 &\sigma_k \\ \sigma_k &0 \\ \end{matrix}\right) \qquad (k=1, 2, 3)\;,
$$
where ${
\mathbb I}=\left( \begin{matrix} 1 & 0 \\ 0 & 1 \end{matrix} \right)$ and $\sigma_k$ are the Pauli matrices:
$$\sigma _1=\left( \begin{matrix} 0 & 1 \\ 1 & 0 \\ \end{matrix} \right),\quad \sigma_2=\left( \begin{matrix} 0 & -i \\ i & 0 \\ \end{matrix}\right),\quad \sigma_3=\left( \begin{matrix} 1 & 0\\ 0 &-1\\ \end{matrix}\right) .$$

The magnetic Dirac operator without the Coulomb potential has essential spectrum $(-\infty, -1] \cup [1, \infty)$ and no eigenvalues in the gap $(-1,1)$. For $\nu\in(0,1)$ the Hamiltonian (\ref{one}) has the same essential spectrum and eigenvalues in the gap. The ground state energy $\lambda_1(\nu, B)$ is defined as the smallest eigenvalue in the gap.

As the field gets large enough, one expects that the ground state energy of the Dirac operator decreases and eventually penetrates the lower continuum. The implication of this for a second quantized model is that electron--positron pair creation comes into the picture {\cite{0627.58040,pickl}}. The intuition for that can be gleaned from the Pauli equation, where the magnetic field tends to lower the energy because of the spin. It is therefore reasonable to define the {\it critical field strength} $B(\nu)$ as the supremum of the {positive} $B$'s for which $\lambda_1(\nu, b)$ is in the gap $(-1,1)$ for all {$b\in (0,B)$}. As a function of $\nu$, $\lambda_1(\nu, B)$ is non-increasing. Hence the function $B(\nu)$ is also non-increasing.

\medskip One of our goals is to give estimates on this critical field as a function of the nuclear charge. Our first result, proved in Section 2 is that this critical field exists and we give some rough estimates in terms of $\nu$: For some $C>0$,
\be\label{Estimation1}
{\frac {4}{5\,\nu^2}} \leq\; B(\nu)\;\leq \;\min \left(\frac{18\,\pi \nu^2}{[3\,\nu^2-2]_+^2}\ ,\ e^{\,C/\nu^2} \right)\,.
\ee
As a corollary we get the noteworthy result that as $\nu \to 1$ the critical field $B(\nu)$ stays strictly positive. This is somewhat remarkable, since in the case without magnetic field the ground state energy as a function of $\nu$ tends to $0$ as $\nu \to 1$ but {\it with an infinite slope}. Thus, one might expect very large variations of the eigenvalue at $\nu=1$ as the magnetic field is turned on, in particular one might naively expect that the ground state energy leaves the gap for small fields $B$. This is not the case. Moreover, since the hydrogenic Hamiltonian ceases to be selfadjoint at $\nu=1$ it is hard to visualize how one might arrive at such estimates using standard perturbation theory.

\medskip Section 3 is devoted to the asymptotics of $B(\nu)$ as $\nu\to 0$. We define the notion of lowest relativistic Landau level which leads to {a} one dimensional effective theory. This effective theory can be analyzed in great detail and allows to calculate the ground state energy $\llc(\nu,B)$ of the magnetic Dirac--Coulomb equation~(\ref{one}) in the lowest relativistic Landau level. It is given by the variational problem
{\[
\llc(\nu,B) := \inf_{f\in C^\infty_0(\R, \C)\setminus\{0\}}\quad \lambdaL[f,\nu,B]\ ,
\]}
where {$\lambda=\lambdaL[f,\nu,B]$} is defined by
\[
\lambda\int_\R |f(z)|^2\,dz= \int_\R\left(\frac{|f'(z)|^2}{1+ \lambda+\nu\,a^B_0(z)}+(1-\nu\,a^B_0(z))\,|f(z)|^2\right)dz\;,
\]
and
$$
{a^B_0(z)}:=B\,\int_0^{+\infty}\frac{s\,e^{-\frac12\,B\,s^2}}{\sqrt{s^2+z^2}}\,ds\;.
$$
The point here is that for $B$ not too small and $\nu$ not too large (the precise bounds are given in Theorem \ref{thmfirstll}),
\[
\llc(\nu+\nu^{3/2}, B) \le \lambda_1(\nu, B) \le \llc(\nu-\nu^{3/2},B) \ .
\]

The one dimensional $\llc(\nu,B)$ problem, although not trivial, is simpler to calculate than the $\lambda_1(\nu, B)$ problem. As a result, in the limit as $\nu \to 0$, this new theory yields the first term in the asymptotics of the logarithm of the critical field. In particular we have the following result,
$$
\lim_{\nu \to 0} \nu \log(B(\nu))\,=\,\pi\ .
$$

\medskip From a methodological point of view, the ground state energy of the Dirac operator is {\it not} given by a minimum problem for the corresponding Rayleigh quotient, but it is a min--max in the sense that one decomposes the whole Hilbert space ${\mathcal H} = {\mathcal H}_1 \oplus {\mathcal H}_2$, maximizes the energy over functions in ${\mathcal H}_2$ and then minimizes over non-zero functions in ${\mathcal H}_1$. While the choice of these Hilbert spaces is not arbitrary, there is some flexibility in choosing them, see \cite{GS,GLS,DES}. For certain choices, the maximization problem can be worked out almost explicitly leading to a new energy functional for which the ground state energy $\lambda_1(\nu, B)$ is the minimum. In this sense, the ground state energy of the Dirac operator appears as a minimum of a well defined functional. Both variational characterizations, the min--max and the min, are of course equivalent, and our approach depends on the interplay between the two.

\medskip Our results are different from the work of \cite{AHS} which considered the non-relativistic hydrogen atom and worked out the asymptotics of the ground state energy as $B \to \infty$ for every $\nu >0$. In our case, however, $\nu$ has to stay in the interval $[0,1)$ in order that the operator can be defined as a selfadjoint operator. Further, the critical field is always finite and we are interested in estimating it as a function of $\nu$. The similarity with \cite{AHS} comes as we let $\nu \to 0$, since then the critical field tends to infinity but the estimates are not the same.

\medskip While the mathematical methods are the main point of this paper, let us make a few additional remarks about its physical motivation. Spontaneous pair creation in strong external fields, although never experimentally confirmed, has been analyzed by Nenciu \cite{N,0627.58040}. In \cite{N} it was conjectured that by adiabatically switching the potential on and off, there is spontaneous pair creation provided some of the eigenvalues emerging from the negative spectrum crossed eigenvalues emerging from the positive spectrum. This conjecture was partly proved in \cite{0627.58040} and \cite{pickl}. Since such a crossover occurs in the Dirac hydrogenic atom with a strong magnetic field, it is natural to try to estimate the strength of the magnetic fields for which this crossing phenomenon occurs. 

Note that the unit in which we measure the magnetic field is huge, about {$4.4\times 10^{9}$} Tesla. Sources of gigantic magnetic fields are neutron stars that can carry magnetic fields of about $10^{9}$ Tesla. Fields of $10^{11}$ Tesla for a neutron star in its gestation and in magnetars are expected, and there is speculation that fields of up to $10^{12}$ Tesla may exist in the interior of a magnetar. There is a considerable literature in this area and an entertaining introduction can be found in~\cite{ScientificAmerican}.

Further, it is expected that near the surface of a neutron star atoms persist up to about $Z=40$. We show that the critical field at $Z=40$ must be larger than $4.1 \times 10^{10}$ Tesla, and preliminary calculations using numerical methods based on Landau levels yield an estimated value of about {$2.5 \times 10^{16}$ Tesla}. Although improvements on these estimates are currently under investigation we believe it is unlikely that they will yield relevant values for the magnetic field strength. For elements with higher $Z$, the values for the critical field are much lower. In the case of Uranium ($Z=92$), they are sandwiched between $7.8\times 10^9$ Tesla and an estimated value (using Landau levels) of $4.6\times 10^{11}$ Tesla.

\medskip Speculations that large magnetic fields facilitate the creation of electron - positron pairs are not new in the physics and astrophysics literature. Clearly, the Dirac operator coupled to a magnetic field but {\it without} electrostatic potential has a gap of $2mc^2$ independent of the magnetic field. It was pointed out in \cite{PhysRev.173.1220,PhysRevLett.21.397} that the anomalous magnetic moment narrows the gap, {\sl i.e.,\/} it decreases the energy needed for pair production. In lowest nontrivial order the anomalous magnetic energy is proportional to the magnetic field which leads indeed to a narrowing of the gap; in fact the gap closes at a field strength of about $4 \times 10^{12}$ Tesla. It was observed in \cite{PhysRev.187.2275}, however, that the anomalous magnetic energy depends in a non linear fashion on the external field. Further it is shown that even at field strengths of $10^{12}$ Tesla the gap narrows only a tiny bit, irrelevant for pair production. For a review of these issues the reader may consult \cite{duncan}. Our contribution is to take into account simultaneously the magnetic field and the Coulomb singularity, in which case no explicit or simple calculations are possible.

Of course our analysis only deals with a single electron and a fixed nucleus. A description of the non-relativistic many electron atom under such extreme situations has been given in \cite{MR1163415,MR1272387,MR1266071}. The authors study various limits as the nuclear charge and the magnetic field strength gets large and determine the shape of the atom in these limits. In non-relativistic physics the natural scale for the magnetic field is {$\alpha^2 \frac{m^2 c^2}{|q| \hbar} \approx 2.4\times 10^5$ Tesla}, much smaller than the ones under considerations in our paper.
 
As we have mentioned before, for small $Z$ the critical magnetic field is of the order of $e^\frac{\pi}{Z\alpha}$ and hence non relativistic physics is sufficient to explain the rough shapes and sizes of atoms even at very high field strengths. It may very well be, however, that for heavy elements and very large fields qualitatively new effects appear that cannot be understood on the basis of non relativistic physics alone. Should such effects occur, then it could make sense to treat the many body relativistic electron problem using the Dirac - Fock approximation.

\section{Ground state and critical magnetic field}\label{Sec:Estim}

{In this section, we set some notations, establish basic properties and prove Estimate~\eqn{Estimation1} on the critical magnetic field.}

\subsection{Min--max characterization of the ground state energy}

The eigenvalue equation for the Hamiltonian (\ref{one})
\be \label{eigenvalue}
H_B\,\psi- \frac{\nu}{|x|}\,\psi = \lambda\,\psi
\ee
is an equation for four complex functions. It is convenient to split $\psi$ as
\[
\psi=\binom{\,\phi\,}\chi
\]
where $\phi, \chi \in L^2(\R^3 ; \C^2)$ are the {{\sl upper and lower components.\/}} Written in terms of {$\phi$ and $\chi$,} (\ref{eigenvalue}) is given by
\begin{eqnarray}
P_B\chi + \phi - \frac{\nu}{|x|}\,\phi = \lambda \phi \label{1st} \ , \\
P_B\phi -\chi - \frac{\nu}{|x|}\,\chi = \lambda \chi \label{2nd} \ .
\end{eqnarray}
Here $P_B$ denotes the operator
{
\[
P_B:=-\,i\,\sigma \cdot (\nabla-i\,{ A}_B(x))\;,
\]
}
where
\[
{ A}_B(x):=\frac B2\left(\begin{array}{c} -x_2\\ x_1\\ 0\end{array}\right)
\]
is the magnetic potential associated with the constant magnetic field
\[
{ {\mathbf B}}(x):=\left(\begin{array}{c} 0\\ 0\\ B\end{array}\right)\,.
\]

Using (\ref{2nd}) we can eliminate the lower component $\chi$ in (\ref{1st}). Taking then the inner product with $\phi$ we get
\be
J[\phi,\lambda,\nu,B]=0\ , \label{jayequation}
\ee
where
\[
J[\phi,\lambda,\nu,B]:= \int_{\R^3}\left(\frac{|P_B \phi|^2} {1+\lambda+\frac\nu{|x|}}+(1-\lambda)|\phi|^2-\frac\nu{|x|}\,|\phi|^2\right)\,d^3x\;. \label{jay}
\]
Thus, we see that that for any eigenvalue $\lambda \in (-1,1)$ of (\ref{one}) the corresponding eigenvector leads to a solution of (\ref{jayequation}).

{Reciprocally, the functional $J$ can be used to characterize the eigenvalues. For this purpose,} a few definitions are {useful.} The functional $J[\phi,\lambda,\nu,B]$ is defined for
any $B\in\R^+$, $\nu\in (0,1)$, $\lambda\geq -1$ and $\phi\in C^\infty_0(\R^3, \C^2)$. Further, in order that (\ref{jayequation}) makes sense we introduce the set
\begin{eqnarray*}
{\cal A}(\nu,B):=\{\phi\in C^\infty_0(\R^3)\,:&&\|\phi\|_{L^2(\R^3)}=1,\\ &&\qquad\lambda\mapsto J[\phi,\lambda,\nu,B]\;\mbox{changes sign in}\;(-1,+\infty)\} \,.
\end{eqnarray*}
Note that this set might be {\sl a priori\/} empty. Finally, since the function $J$ is decreasing in $\lambda$, we define $\lambda=\lambda[\phi,\nu,B]$ to be either the unique solution to
\[
J[\phi,\lambda,\nu,B]=0\; \hbox{if}\; \phi\in{\cal A}(\nu,B)\ ,
\]
or $\lambda[\phi,\nu,B]=-1$ if $J[\phi,-1,\nu,B]\leq 0$.
\begin{theorem}\label{min-max} Let $B\in\R^+$ and $\nu\in (0,1)$. If
\[
-1<\inf_{\phi\in C^\infty_0(\R^3, \C^2)}\lambda[\phi,\nu,B] <1\ ,
\]
this infimum is achieved and
\[
\lambda_1(\nu,B):=\inf_{\phi\in{\cal A}(\nu,B)}\lambda[\phi,\nu,B]
\]
is the lowest eigenvalue of $H_B-\frac{\nu}{|x|}$ in the gap of {its} continuous spectrum, $(-1,1)$. \end{theorem}
\proof This proposition is a consequence of Theorem 3.1 in \cite{DES}. The essential assumptions of this theorem are:
\\
i) The selfadjointness of $H_B-\nu\,|\cdot|^{-1}$ which is proved in the appendix. It is crucial here that $0 < \nu < 1$.
\\
ii) {The existence of} a direct decomposition of $L^2(\R^3;\C^4)$ as the sum of two subspaces ${\mathcal H}_1 \oplus {\mathcal H}_2$ such that
\begin{eqnarray*}
a_2 := &&\sup_{x \in {\mathcal H}_2} \frac {(x, (H_B-\nu\,|\cdot|^{-1})\,x)}{(x,x)} \\
&&< c_1 := \kern -2pt\inf_{0 \not= x \in {\mathcal H}_1} \ \sup_{y \in {\mathcal H}_2} \frac {(x+y, (H_B-\nu\,|\cdot|^{-1})\,(x+y))}{\| x+y \|^2} \ .
\end{eqnarray*}
Set $b:= \inf \sigma_{ess}(H_B-\nu\,|\cdot|^{-1}) \cap (a_2, +\infty)$. If $c_1 < b$ then $c_1$ is the lowest eigenvalue of $H_B-\nu\,|\cdot|^{-1}$ in the interval $(a_2, b)$. {In the present case} we choose the decomposition
$$
\psi = \binom\phi\chi = \binom\phi{0}+\binom 0\chi
$$
{based on the upper and lower components of the four components spinor $\psi$.} It is easy to see that $a_2=-1$. Furthermore the essential spectrum of $H_B-\nu\,|\cdot|^{-1}$ is $(-\infty, -1] \cup [1,+\infty)$ independently of $B$, see \cite{T}. Hence $b=1$. It remains to calculate the {supremum} in the definition of $c_1$ as a function of $x = \binom\phi{0}$. Note that the Rayleigh quotient in the definition of $c_1$ is strictly concave in $y=\binom 0\chi$. Therefore the supremum is uniquely achieved by
\[\label{maxx}
\chi[\phi]= \left(1+\lambda[\phi,\nu,B]+\frac{\nu}{|x|}\right)^{-1}P_B\phi
\]
and its value is $\lambda[\phi, \nu, B]$, that is,
\[\label{firstminmax}
\lambda[\phi,\nu,B] = {\sup_{\chi\in C^\infty_0(\R^3, \C^2),\,\psi=\binom\phi\chi}}\;\frac{\Big(\big( H_B-\frac{\nu}{|x|}\big)\,\psi, \psi\Big)}{(\psi, \psi)}\ .
\]
\finprf
\noindent{\bf Remarks:} 1) {\sl The eigenvalue $\lambda_1(\nu, B)$ can be characterized either as the minimum of the functional $\lambda[\phi,\nu,B]$ or as a min-max level of $H_B-\nu\,|\cdot|^{-1}$. {Both characterizations will be} useful in the sequel of this paper.}

\smallskip\noindent 2) {\sl Under the assumptions of Theorem \ref{min-max}, we have
\[
J[\phi, \lambda,\nu, B]\geq 0\quad \forall \phi\in C^\infty_0(\R^3, \C^2)\label{newrm}
\]
for any $\lambda \le \lambda_1(\nu, B)$. {The eigenvalue $\lambda_1(\nu, B)$ can therefore be interpreted as the best constant in the above inequality.}}

\smallskip\noindent 3) {\sl When $\lambda_1(\nu,B)$ is equal to $-1$, it belongs to the continuous spectrum and it is not necessarily an eigenvalue of $H_B-\nu\,|\cdot|^{-1}$.}

\subsection{Basic properties of the ground state energy}

\begin{proposition}\label{Prop:MonotonicityNu} For all $B\geq 0$, the function $\nu\mapsto \lambda_1(\nu,B)$ is monotone nonincreasing on $(0,1)$. \end{proposition}
{The proof is left to the reader. It is a consequence of the definition of $J[\phi,\lambda,\nu,B]$. 
\begin{proposition}\label{cont} For all $B\geq 0$, the function $\nu\to \lambda_1(\nu,B)$ is continuous in the interval $\left(0,1\right)$ as long as $ \lambda_1 (\nu,B)\in (-1,1)$. \end{proposition}
\proof By Theorem \ref{min-max}, if $ \lambda_1 (\nu,B)\in (-1,1)$ there exists a function $\phi_\nu$ such that $J[\phi_\nu, \lambda_1(\nu,B), \nu, B]=0$. For any sequence $\{\nu_n\}_{_n}$ converging to~$\nu$, the upper semi-continuity of $\nu\to \lambda_1(\nu,B)$ holds:
\[\label{uppercont}
\limsup_{n\to +\infty} \lambda_1(\nu_n, B)\leq \limsup_{n\to +\infty} \lambda[\phi_\nu,\nu_n, B]=\lambda_1(\nu, B)\ .
\]
{If $\nu_n\leq\nu$, then $\lambda_n:=\lambda_1(\nu_n,B)\geq\lambda_1(\nu,B)$ and $\{\lambda_n\}_{_n}$ converges to $\lambda$. Consider therefore} a $\{\nu_n\}_{_n}$ converging to $\nu$ from above. {We have to face two~cases:\\
\underline{\sl First case:\/} $\lambda_n>-1$ for all $n\in\N$. Since $J[\phi_\nu, \lambda_1(\nu,B), \nu_n, B]\leq 0$, we know that $\lambda_n\leq\lambda_1(\nu,B)$.} Consider the corresponding eigenfunctions $\psi_n$, such that $J[\phi_n, \lambda_n, \nu_n, B]=0$, where $\phi_n$ denotes the upper component of $\psi_n$ {and assume that $\|\phi_n\|_{L^2(\R^3)}=1$. By Theorem~\ref{min-max}, we have}
\[
\int_{\R^3}\left(\frac{|P_B \phi_n|^2} {1+\lambda_n+\frac{\nu_n}{|x|}} +(1-\lambda_n)|\phi_n|^2\right)\,d^3x =\int_{\R^3} \frac{\nu_n}{|x|}\,|\phi_n|^2 \,d^3x\;.
\]
{Assume that $\Lambda:=\liminf_{n\to +\infty} \lambda_1(\nu_n,B)<\lambda_1(\nu,B)$. Up to the extraction of a subsequence, assume further that $\{\lambda_1(\nu_n,B)\}_{n\in\N}$ converges to some value in $(-1,\lambda_1(\nu,B))$ and choose $\tilde\lambda\in(\Lambda,\lambda_1(\nu,B))$. For $n$ large enough, $\lambda_1(\nu_n,B)<\tilde\lambda$ and 
\be\label{nunb}
\int_{\R^3}\left(\frac{|P_B \phi_n|^2} {1+\tilde\lambda+\frac{\nu_n}{|x|}}+(1-\tilde\lambda)|\phi_n|^2\right)\,d^3x \leq\int_{\R^3} \frac{\nu_n}{|x|}\,|\phi_n|^2 \;.
\ee
} \par\smallskip\noindent
{\underline{\sl Second case:\/}} $\lambda_1(\nu', B)=-1$ for all $ \nu'>\nu$. We choose $\tilde\lambda\in (-1, \lambda_1(\nu, B))$ and find a $\{\phi_n\}_n$ such that $\|\phi_n\|_{L^2(\R^3)}=1$ and $ \,J[\phi_n, \tilde\lambda,\nu_n, B]\leq 0$ for $n$ large: \eqn{nunb} also holds. 

\par\medskip Using the monotonicity of the $\{\nu_n\}_n$, {which implies the monotonicity of the $\{\lambda_n\}_n$ by Proposition \ref{Prop:MonotonicityNu}, and the fact that in both cases, $\nu\in(0,1)$ and $\tilde\lambda\in(-1,1)$, we deduce from} (\ref{nunb}) a uniform bound for the functions~$\phi_n$:
\be\label{zzz}
\sup_n\int_{\R^3}\left(|x|\,|P_B \phi_n|^2+\frac{| \phi_n|^2}{|x|}\right)\,d^3x<+\infty\ .
\ee
{The proof goes as follows. It is sufficient to prove that $\int_{\R^3}|x|^{-1}\,|\phi_n|^2\,d^3x$ is uniformly bounded. Let $\chi$ be a smooth truncation function such that $\chi(r)\equiv 1$ if $r\in [0,1)$, $\chi(r)\equiv 0$ if $r>2$, and $0\leq \chi\leq 1$. Since
\begin{eqnarray*}
\int_{\R^3}\frac{|\phi_n|^2}{|x|}\,d^3x&\leq& \int_{\R^3}\frac{|\tilde \phi_n|^2}{|x|}\,d^3x+\frac 1R\int_{\R^3}|\phi_n|^2\,\left(1-\chi^2\left(\frac{|x|}R\right)\right)\,d^3x\\
&\leq& \int_{\R^3}\frac{|\tilde \phi_n|^2}{|x|}\,d^3x+\frac 1R
\end{eqnarray*}
with $\tilde\phi_n=\chi(|R^{-1}\cdot|)\,\phi_n$, it is therefore sufficient to prove that $\int_{\R^3}|x|^{-1}\,|\tilde\phi_n|^2\,d^3x$ is uniformly bounded, for some $R>0$, eventually small. Using the estimate
\[
a^2\geq \frac{(a+b)^2}{1+\varepsilon}-\frac{b^2}\varepsilon\ ,
\]
we get the following lower bound
\begin{eqnarray*}
\int_{\R^3}\frac{|P_B\phi_n|^2}{1+\tilde\lambda+\frac{\nu_n}{|x|}}\,d^3x&\geq&
\int_{\R^3}\frac{|P_B\phi_n|^2\,\chi^2}{1+\tilde\lambda+\frac{\nu_n}{|x|}}\,d^3x\\ &\geq& \int_{\R^3}\frac{|P_B\tilde\phi_n|^2}{(1+\varepsilon)\big(1+\tilde\lambda+\frac{\nu_n}{|x|}\big)}\,d^3x-\frac C\varepsilon\,\|\tilde\phi_n\|_{L^2(\R^3)}^2
\end{eqnarray*}
}{where $C$ is a constant which depends on $\|\chi'\|_{L^\infty(\R^+)}^2$, $B$ and $R$. Next, with the same type of arguments, we can write
\[
|P_B\tilde\phi_n|^2\geq\frac{|\sigma\cdot\nabla\tilde\phi_n|^2}{1+\varepsilon}-\frac 1\varepsilon\,{\Big|\,B\,|x|\,\tilde\phi_n\,\Big|^2}\geq\frac{|\sigma\cdot\nabla\tilde\phi_n|^2}{1+\varepsilon}-\frac {B^2\,R^2}\varepsilon\,|\tilde\phi_n|^2\;.
\]
Collecting these estimates, this gives
\[
\int_{\R^3}\frac{|\sigma\cdot\nabla\tilde\phi_n|^2}{(1+\varepsilon)^2\big(1+\tilde\lambda+\frac{\nu_n}{|x|}\big)}\,d^3x\leq C(\varepsilon,R,\chi)+\int_{\R^3} \frac{\nu_n}{|x|}\,|\tilde\phi_n|^2 \,d^3x\;.
\]
Because $\tilde\phi_n$ has a compact support in the ball of radius $2\,R$, if $\delta=\delta(R)>(1+\tilde\lambda)R/\nu_n$ at least for $n$ large enough, then 
\[
\frac 1{1+\tilde\lambda+\frac{\nu_n}{|x|}}\geq \frac{|x|}{\nu_n\,(1+\delta)}\quad\forall\;x\in B(0,R)
\]
so that
\[
\frac1{(1+\varepsilon)^2\,\nu_n\,(1+\delta)}\int_{\R^3}|x|\,|\sigma\cdot\nabla\tilde\phi_n|^2\,d^3x\leq C(\varepsilon,R,\chi)+\int_{\R^3} \frac{\nu_n}{|x|}\,|\tilde\phi_n|^2 \,d^3x\;.
\]
On the other hand, according to \cite{DES,DELV},
\[
\int_{\R^3}|x|\,|\sigma\cdot\nabla\tilde\phi_n|^2\,d^3x\geq \int_{\R^3} \frac 1{|x|}\,|\tilde\phi_n|^2 \,d^3x\ .
\]
This provides a uniform upper bound on $\int_{\R^3}|x|^{-1}\,|\tilde\phi_n|^2 \,d^3x$ if $\varepsilon$ and $\delta$ are chosen small enough in order that 
\[
\frac1{(1+\varepsilon)^2\,(1+\delta)}>\nu_n^2
\]
for $n$ large. This can always be done since $\nu_n$ converges to $\nu\in (0,1)$ and $\delta(R)$ can be taken as small as desired for $R>0$ sufficiently small. This concludes the proof of \eqn{zzz}.}

\medskip{Summarizing,} we obtain that
$$
J\left[\phi_n, {\frac 12(\tilde\lambda+\lambda_1(\nu,B))}, \nu, B\right] \leq 0\ ,
$$
for {$n$ large enough: hence} $\lambda [\phi_n, \nu, B]\leq\frac 12{(\tilde\lambda+\lambda_1(\nu,B))}<\lambda_1(\nu, B)$, a contra\-diction. \finprf

Consider now the effect of a scaling on $J$.}
\begin{lemma} Let $B\geq 0$, $\lambda\geq -1$, $\theta>0$ and $\phi_\theta(x) :=\theta^{3/2}\phi(\theta\,x)$ for any $x\in\R^3$. Then
\[
\nabla_{{ A}_{\theta^2B}}\phi_\theta(x)=\theta^{5/2}\left[\nabla\phi(\theta\,x)-i\,{ A}_B(\theta\,x)\,\phi (\theta\,x) \right] \ ,
\]
and for any $a \in \R$, $\nu \in (0,1)$,
\be\label{Eqn:ScaledFunctional}
J[\phi_\theta,\lambda,\theta^a\nu,\theta^2B]= \int_{\R^3}\left(\theta^2\frac{|P_B \phi|^2} {1+\lambda+\frac{\theta^{a+1}\nu}{|x|}} +(1-\lambda)|\phi|^2-\frac{\theta^{a+1}\nu}{|x|}\,|\phi|^2\right)\,d^3x\;.
\ee \end{lemma}
Using this scaling, we prove some properties enjoyed by the function $\lambda_1(\nu,B)$.
Take $\lambda=\lambda[\phi_\theta,\nu,\theta^2B]>-1$, $ \theta>1$ and $a=0$ in \eqn{Eqn:ScaledFunctional}. With $1+\lambda= \theta(1+\mu)$,
\begin{eqnarray*}
0&=&J[\phi_\theta,\lambda,\nu,\theta^2B]\\
&=&\theta\int_{\R^3}\left(\frac {|P_B\phi|^2} {1+\mu+\frac\nu{|x|}}+(1-\mu)|\phi|^2-\frac\nu{|x|}\,| \phi|^2\right)d^3x+2(1-\theta)\int_{\R^3}|\phi|^2\;d^3x\ .
\end{eqnarray*}
Assuming that $\|\phi\|_{L^2(\R^3)}=1$, we get
\[
J[\phi,\mu,\nu,B]=2\,\frac{\theta-1}\theta>0\,
\]
and thus,
\[
\lambda[\phi,\nu,B]> \mu=\frac\lambda\theta-\frac{\theta-1}\theta\;.
\]
On the other hand, $\frac \partial{\partial\mu}J[\phi,\mu,\nu,B]\leq -1$, so that
\[
\lambda[\phi,\nu,B]\leq \mu+J[\phi,\mu,\nu,B]=\frac\lambda\theta+\frac {\theta-1}\theta\;.
\]
Summarizing, we have the estimate
\[
\frac{\lambda[\phi_\theta,\nu,\theta^2B]}\theta-\frac{\theta-1}\theta \leq \lambda[\phi,\nu,B] \leq \frac{\lambda[\phi_\theta,\nu, \theta^2B]}\theta+\frac{\theta-1}\theta\quad\forall\; \theta>1\;.
\]
The above estimate, which holds provided $\lambda[\phi_\theta,\nu,\theta^2B]>-1$ is equivalent to
\be\label{Estim-Page5}
\frac{\lambda[\phi,\nu,\theta^2B]}\theta-\frac{\theta-1}\theta\leq \lambda[\phi_{1/\theta},\nu,B] \leq \frac{\lambda[\phi,\nu,\theta^2B]} \theta+\frac{\theta-1}\theta\quad\forall\; \theta>1
\ee
under the condition $\lambda[\phi,\nu,\theta^2B]>-1$. {As a consequence, we have the following result.}
\begin{proposition}\label{Prop:Lipschitz} For all $\nu\in (0,1)$, the function $B\mapsto \lambda_1(\nu,B)$ is continuous as long as it takes its values in $(-1,+\infty)$. Moreover
\be\label{restricttheta}
\frac{\lambda_1(\nu,\theta^2B)}\theta-\frac{\theta-1}\theta\leq \lambda_1(\nu,B) \leq \frac{\lambda_1(\nu,\theta^2B)}\theta+\frac {\theta-1}\theta\ ,
\ee
if $\lambda_1(\nu, B) \in (-1, +\infty)$ and $\theta\in\left(1, \frac 2{1-\lambda_1(\nu,B)}\right)$. As a consequence,
$B\mapsto \lambda_1(\nu,B)$ is Lipschitz continuous for any $\nu \in (0,1)$ and $B>0$ such that $\lambda_1(\nu,B)\in(-1,+\infty)$:
\[
\frac{\lambda_1-1}{2B} \leq\frac{\partial\lambda_1}{\partial B}\leq \frac{\lambda_1+1}{2B}\;. \]\end{proposition}
\proof Choose $a\in (-1,\lambda_1(\nu,B))$ and take any $\phi\in C^ \infty_0(\R^3, \C^2)$. Since
$$
\frac\partial{\partial\lambda}J[\phi, \lambda,\nu,B]\leq -1\ ,
$$
an integration on the interval $[a,\lambda [\phi,\nu,B]]$ shows that
\[
-J[\phi,a,\nu,B]=\Big[J[\phi,\lambda,\nu,B]\Big]_{\lambda=a}^{\lambda= \lambda[\phi,\nu,B]}\leq -\lambda[\phi,\nu,B]+a
\]
where the first {equality} holds by definition of $\lambda[\phi,\nu,B]$, {\sl i.e.\/} $J[\phi,\lambda[\phi,\nu,B],\nu,B]=0$. As a consequence,
\[
J[\phi,a,\nu,B]\geq\lambda[\phi,\nu,B]-a>0\;.
\]
The function $\theta\mapsto J[\phi,a,\nu,\theta^2B]$ is continuous, so only two cases are possible:\\
\underline{\sl First case:\/}
\[
J[\phi,a,\nu,\theta^2B]\geq 0\quad\forall\;\theta>1\;,
\]
\underline{\sl Second case:\/} there exists a constant $\bar\theta=\bar\theta(a,\phi)>1$ such that \begin{description}
\item{(i)} $J[\phi,a,\nu,\theta^2B]>0$ for any {$\theta\in (1,\bar \theta)$,}
\item{(ii)} $J[\phi,a,\nu,{\bar\theta}^2B]=0$ or, equivalently, $ \lambda[\phi,\nu,{\bar\theta}^2B]=a$.
\end{description}
In the second case, by (i), we know that $\lambda[\phi,\nu,\theta^2B] > a>-1$ for any {$\theta\in (1,\bar\theta)$} and so \eqn{Estim-Page5} applies:
{\[
\theta\,\lambda_1(\nu,B)\leq \theta\,\lambda[\phi_{1/\theta},\nu,B]\leq\lambda[\phi,\nu,\theta^2B]+\theta-1\;.
\]}
In the limit case $\theta=\bar\theta$, we get
\[
\bar\theta\,\lambda_1(\nu,B)+1-\bar\theta\leq\lambda[\phi,\nu,{\bar \theta}^2B]=a\;,
\]
which gives the estimate
\[
\bar\theta\geq\frac{1-a}{1-\lambda_1(\nu,B)}=:\theta^*(a)\;.
\]
Thus the inequality
\[
\theta\,\lambda_1(\nu,B)\leq \theta\,\lambda[\phi_{1/\theta},\nu,B]\leq\lambda[\phi,\nu,\theta^2B]+\theta-1
\]
holds for any $\theta\in [0,\theta^*(a)]$ and for any $\phi\in C^ \infty_0(\R^3, \C^2)$, which proves the r.h.s. inequality in (\ref {restricttheta}) by letting $a\to -1$:
\[
\lim_{a\to -1}\theta^*(a)=\frac 2{1-\lambda_1(\nu,B)}\;.
\]
The l.h.s. inequality is obtained in the same manner.
\finprf

With the appropriate test functions one can prove that $\lambda_1(\nu, B)$ is always below~$-1$ for $B$ large. We recall that $\lambda_1(\nu,B)=-1$ means that if $J[\phi,-1,\nu,B]\leq 0$ for any $\phi\in{\cal A}(\nu,B)$. Let us give some details.
\begin{proposition}\label{Prop:testfunction} Let $\nu\in (0,1)$. Then {for $B$ large enough, $\lambda_1(\nu,B)\leq 0$ and there exists $B^*>0$ such that $\lambda_1(\nu,B)=-1$ for any $B\geq B^*$.} \end{proposition}
\proof Let us consider $B>0$ and the trial function
$$
\psi=\binom\phi{0}\ ,
$$
where
$$
\phi = \sqrt{\frac B{2\,\pi}}\,e^{-\frac{B}{4}(|x_1|^2+|x_2|^2)} \,\binom{f (x_3)}{0} \ $$
and $f\in C^\infty_0(\R, \R)$ is such that $f\equiv 1$ for $| x|\leq \delta$, $\delta$ small but fixed, and $\|f\|_{_{L^2 (\R)}}=1$. Note that $\phi\in Ker(P_B+i\sigma_3\,\partial_{x_3})$ and so, 
$$
P_B \phi = -i\,\sigma_3\,\partial_{x_3}\,\phi\; ,
$$
Moreover, the state $\phi$ is normalized in $L^2(\R^3)$. With $r=|x|$, we can define
\be\label{GB}
G_B[\phi]:=\int_{\R^3}\left(\frac r{\nu}\,|P_B \phi|^2- \frac{\nu} {r}\,|\phi|^2\right)d^3x
\ee
{and compute}
\begin{eqnarray*}
G_B[\phi]&=&B\,\iint_{\R\times\R^+} \left(\frac{\sqrt{s^2+|x_3|^2}}{\nu} \,|f'(x_3)|^2-\frac{\nu\,|f(x_3)|^2}{\sqrt{s^2+|x_3|^2}} \right)s\,e^ {-B\,s^2/2}\,ds\,dx_3
\\&=&\iint_{\R\times\R^+} \left(\frac{\sqrt{\scriptstyle B^{-1}s^2+|x_3|^2}}{\nu} \,|f'(x_3)|^2-\frac{\nu\,|f(x_3)|^2}{\sqrt{\scriptstyle B^{-1}s^2+|x_3|^2}} \right)s\,e^{-s^2/2}\,ds\,dx_3\ .
\end{eqnarray*}
{Using $\sqrt{B^{-1}s^2+|x_3|^2}\leq B^{-1/2}\,s+|x_3|$ and
\begin{eqnarray*}
\iint_{\R\times\R^+}\frac{|f(x_3)|^2}{\sqrt{B^{-1}s^2+|x_3|^2}}\,s\,e^{-s^2/2}\,ds&\geq&
\int_0^1\frac{ds}{\sqrt e}\int_0^\delta\frac s{B^{-1/2}s+x_3}\,dx_3\\
&&\qquad\geq\frac 1{4\sqrt e}\,\log(\delta^2B)\ ,
\end{eqnarray*}
for $B\geq 1$, we can therefore bound $G_B[\phi]$ from above by
\[
G_B[\phi]\leq\frac{C_1}{\nu}+C_2\,\nu- C_3\,\nu\log B\ .
\]
where $C_i$, $i=1$, $2$, $3$, are positive constants which depend only on $f$. For $B\geq 1$ large enough,} $G_B[\phi]+2\,\| \phi \|^2 \leq 0$ and $\lambda_1(\nu,B)= -1$, since in this case $J[\phi,-1,\nu,B] \leq 0$. \finprf

\subsection{The critical magnetic field}

Proposition \ref{Prop:testfunction} motivates the following definition.
\begin{definition}\label{Def:CriticalMagField} Let $\nu\in (0,1)$. We define the {\rm critical} magnetic field as
\[
B(\nu):=\inf\left\{B>0\;:\;\lim_{b\nearrow B}\lambda_1(\nu, b)=-1\right\}\;.
\]\end{definition}

\begin{corollary} For all $\nu\in (0,1)$, $\lambda_1(\nu,B)<1$ {for any $B\in(0,B(\nu))$} . \end{corollary}
\proof For $B= 0$ we have $\lambda_1(\nu, 0)=\sqrt{1-\nu^2}<1$. {Given $B>0$, small, by continuity of $B\mapsto\lambda_1(\nu,B)$ we know that $\lambda_1(\nu,B)\in (0,1)$. Let us consider~$\theta\in(1,\sqrt{B(\nu)/B})$ such that $-1<\lambda_1(\nu, \theta^2B)\leq 0$.} This is made possible by Propositions \ref{Prop:Lipschitz} and \ref{Prop:testfunction}. Then, by Proposition \ref{Prop:Lipschitz},
$$ \lambda_1(\nu, B)\leq \frac{\theta -1}{\theta}<1\ .$$ \finprf
The computations {of Proposition \ref{Prop:testfunction}} show the existence of a constant $C_3>0$ such that $B(\nu)\leq e^{\,C_3/\nu^2}$, for all $\nu\in (0,1)$. This estimate can be made more precise for any $\nu$ not too small:
\begin{theorem}\label{HHH} For all {$\nu\in(0,1)$, there exists a constant $C>0$ such that}
\[\label{firstineq}
{\frac {4}{5\,\nu^2}} \leq\; B(\nu)\;\leq \;\min \left(\frac{18\,\pi\,\nu^2}{[3\,\nu^2-2]_+^2}\ ,\ e^{\,C/\nu^2} \right)\ .
\]
\end{theorem}
The proof of this theorem uses Proposition \ref{Prop:testfunction}. Otherwise it is splitted in two partial results stated in Propositions \ref{Prop:testfunction2} and~\ref{Prop:lowerboundatone}. 

Notice that there is a big gap between these lower and upper estimates when $\nu$ is small. To try to better understand this problem, in the next section, we will analyze the case when $B$ is large, proving that the $3d$ definition of $B(\nu)$ is actually asymptotically equivalent to a $1d$ problem related to the lowest relativistic Landau level. More precisely we will prove that when $\nu$ is small, $B$ is not too small and $\lambda_1(\nu, B)\in (-1,1)$, the eigenvalue associated with $\lambda_1(\nu,B)\in (-1,1)$ is {\sl almost equal\/} to the corresponding eigenvalue in the lowest relativistic Landau level class of functions, see Theorem~\ref{thmfirstll}. We will also establish that $B(\nu)$ behaves in the limit $\nu\to 0$ like the upper bound in Theorem~\ref{HHH} and obtain the corresponding value of $C$, see Theorem~\ref{CC}.

Our first partial result is the following
\begin{proposition} \label{Prop:testfunction2} For any $\nu\in (\sqrt{2/3},1)$, $\sqrt{B(\nu)}\leq \frac{3\,\sqrt{2\,\pi}\,\nu}{3\,\nu^2-2}$. \end{proposition}
\proof Consider the trial function $\psi=\binom\phi{0}$ where
$$
\phi = \left(\frac B{2 \pi}\right)^{3/4}e^{-B\,|x|^2/4}\,\binom 10\ ,
$$
is like the one chosen in the proof of Proposition \ref{Prop:testfunction}, {with $f(x_3)=\left(\frac B{2 \pi}\right)^{1/4}e^{-B x_3^3/4}$}. Here, with the notation $r:=|x|$, we find
\begin{eqnarray*}
&&\hspace*{-12pt}G_B[\phi]=B^2\int \frac{r\,|x_3|^2}{4\,\nu}\,|\phi|^2 \,d^3x - \int \frac{\nu}{r}\,|\phi|^2 \,d^3x
\\&=&
(2 \pi)^{-\frac 32}\,B^{\frac 12} \left[\frac{\pi}{2\,\nu} \int_0^\infty e^{-r^2/2}\,r^5\,dr \int_0^\pi\cos^2\theta\,\sin\theta\,d \theta - 4\,\pi\,\nu \int_0^\infty e^{-r^2/2}\,r\,dr\right]
\\&=&
(2 \pi)^{-\frac 32}\,B^{\frac 12}\left[ \frac{\pi}{{3\,\nu}}\int_0^\infty e^{-r^2/2} r^5 dr - 4\,\pi\,\nu \int_0^\infty e^{-r^2/2}\,r\,dr\right]
\\&=&
(2\pi)^{-\frac 32}\,B^{\frac 12}\left[\frac{8\,\pi}{{3\,\nu}}-4\,\pi\,\nu\right]\ .
\end{eqnarray*}
If $\nu^2\in (2/3, 1)$ and $\sqrt{B}\geq \frac{3\,\sqrt{2\,\pi}\,\nu}{3\,\nu^2-2}$, then 
$G_B[\phi]\leq -2=-||\phi||_2$ and so $\lambda_1(\nu, B)=-1$, which proves the Proposition. \finprf

Proposition \ref{Prop:testfunction} shows that for $\nu>\sqrt{2/3}$ and $B$ large, $\lambda_1(\nu, B)$ possibly ceases to be an eigenvalue of the operator $H_B-\nu\,|\cdot|^{-1}$. This can be interpreted by saying that for strong magnetic fields, the Coulomb potential does not stabilize the electron. At some level, electron-positron pairs could appear and then Quantum Field Theory (or QED) becomes unavoidable for a correct description of the electron dynamics, see \cite{pickl}.

\begin{proposition}\label{Prop:ScalingBNu} For given $\nu\in (0, 1)$ and $B>0$, the function {$\theta\mapsto\lambda_1(\theta^{-1}\nu,\theta^2B)$ is monotone nondecreasing as long~as it takes its values in $(-1,1)$ and $\nu/\theta\in(0,1)$.}\end{proposition}
\proof Take $a=-1$ in \eqn{Eqn:ScaledFunctional}:
\[
J\left[\phi_\theta,\lambda, \theta^{-1}\nu,\theta^2B\right]=(\theta^2-1) \int_{\R^3}\frac{|P_B\phi|^2} {1+\lambda+\frac\nu{|x|}}\;d^3x+J[\phi,\lambda,\nu,B]
\]
so that for $\theta<1$,
\[
J\left[\phi_\theta,\lambda,\theta^{-1}\nu,\theta^2B\right]\leq J[\phi,\lambda,\nu,B]
\]
at least for $1-\theta>0$, small, so that $\phi_\theta\in {\cal A}\left(\theta^{-1}\nu,\theta^2B\right)$ for $\phi={\phi_\nu}$ such that $
\lambda[{\phi_\nu},\nu,B] =\lambda_1(\nu,B)$. This proves that
\be\label{Ineq:theta}
\lambda_1\left(\theta^{-1}\nu,\theta^2B\right) \leq\lambda_1(\nu,B)
\ee
for $1-\theta>0$, small. By continuation, the property holds as long as the {assumptions of Proposition~\ref{Prop:ScalingBNu} are satisfied. The case $\theta>1$ follows by multiplying $\theta\,\nu$ and $\theta^{-2}B$ by respectively $\theta^{-1}$ and $\theta^2$.}\finprf
\begin{corollary}\label{Cor:LowerEstimate} There exists a positive constant $\Lambda$ such that
\[
B(\nu)\geq \frac\Lambda{\nu^2}\quad\mbox{as}\;\nu\searrow 0\;.
\] \end{corollary} 
\proof Let $(\nu_0,B_0)$ be such that $B(\nu_0)>B_0$, {\sl i.e.,\/} $\lambda_1(\nu_0,B_0)>-1$ and take $\nu\in (0,\nu_0)$, $\theta=\nu/\nu_0\in (0,1)$, $B=\theta^{-2}B_0$ in \eqn{Ineq:theta}:
\[
-1<\lambda_1(\nu_0,B_0)\leq\lambda_1(\nu,B)=\lambda_1(\sqrt{B_0/B}\,\nu_0,B)\;.
\]
By Proposition \ref{Prop:MonotonicityNu}, this inequality can be extended to
\[
\lambda_1(\sqrt{B_0/B}\,\nu_0,B)\leq \lambda_1(\nu,B)\quad\forall\;\nu\in (0,\sqrt{B_0/B}\,\nu_0)\;.
\]
{This amounts to say} that
\[
0\leq\nu\leq\sqrt{B_0/B}\,\nu_0\Longrightarrow B(\nu)\geq B_0\,\frac{\nu_0^2}{\nu^2}\;,
\]
which proves the result with $\Lambda=B_0\,\nu_0^2$.\finprf

The constant $\Lambda$ can be made more precise. {The remainder of this section is devoted to the following improvement of Corollary \ref{Cor:LowerEstimate}.}
\begin{proposition} \label{Prop:lowerboundatone} For all $\nu\in (0,1)$,
$$
G_B[\phi]=\int_{\R^3}\frac{|x|}{\nu}\,|P_B\phi|^2\,d^3x-\int_{\R^3}\frac{\nu}{|x|}\,|\phi|^2\,d^3x\ge-\nu\,\sqrt {5B}\int_{\R^3}|\phi|^2\,d^3x\ .
$$
In particular this implies that {$B(\nu)\geq \frac 4{5\,\nu^2}$}. \end{proposition}
\proof Scaling the function $\phi$ according to
$$
\phi_B:=B^{3/4}\,\phi\left(B^{1/2}\,x\right)
$$
preserves the $L^2$ norm, and yields
$$
G_B[\phi_B] = \sqrt B\,G_1[\phi] \ ,
$$
{where $G_B$ has been defined in \eqn{GB}.} Obviously it is sufficient to find a good estimate on the functional $G_1$.

\bigskip {Let us collect some preliminary observations.} Recalling that the angular momentum vector $L$ is given by
$$
L= -i \,\nabla \wedge x \ ,
$$
a simple calculation shows that
\[\label{angmom}
\sigma \cdot \nabla = \left(\sigma \cdot \frac{x}{r}\right)\,\left(\partial_r - \frac{1}{r}\,\sigma \cdot L\right)
\]
{with $r=|x|$ and $\partial_r=\frac xr\cdot\nabla$. We also recall that}
$$
-i\,A_1(x) \cdot \sigma = -\left(\sigma \cdot \frac{x}{r}\right) \,\left( \sigma \cdot q(x)\right) \ \quad\mbox{where}\quad q(x) = \frac{1}{2\,r} \left[\begin{array}{c} -x_3\,x_1 \\ -x_3\,x_2 \\ x_1^2+x_2^2 \end{array}\right]
$$
{and}
$$
i\,P_1 = \sigma \cdot \nabla - i\,A_1(x) \cdot \sigma
$$
{so that we can expand $|P_1 \phi |^2$ as
\begin{eqnarray*}
|P_1 \phi |^2&=&\left|\left(\partial_r -\frac{1}{r}\,\sigma \cdot L -\sigma \cdot q(x)\right)\phi\right|^2\\
&=&|\partial_r \phi|^2 + \frac{1}{r^2}\,|\sigma \cdot L\,\phi|^2 + |q|^2\,|\phi|^2 -\partial_r \langle \phi, \sigma \cdot q\,\phi \rangle + \langle \phi, \sigma \cdot (\partial_r q)\,\phi \rangle\\
&&\qquad+\;\frac 1r\,\Big[-\partial_r \langle \phi, \sigma \cdot L\,\phi \rangle+\langle \sigma \cdot L\,\phi, \sigma \cdot q\,\phi \rangle +\langle \sigma \cdot q\,\phi, \sigma \cdot L\,\phi \rangle \Big] \ .
\end{eqnarray*}
As a last preliminary remark, we notice that $r\,\partial_rq = q$.}

\medskip Since the vector potential grows linearly, we localize the problem near the origin. To this end consider the function
$$
t(r)=\left\{\begin{array}{ll} 1 & \hbox{if} \; r \leq R\ , \\ &\\
{R}/{r} & \hbox{if}\; r \geq R\ .\end{array}\right.
$$
Since $t(r)\leq 1$ and $\nu \le 1$ we get the lower bound
\begin{eqnarray*}
G_1[\phi]&\geq&\int_{\R^3} \frac{t(r)\,r}{\nu}\,|P_1 \phi |^2\,d^3 x - \int_{\R^3} \frac{\nu}{r}\,|\phi|^2\,d^3 x
\\&\geq& \nu \left(\int_{\R^3} t(r)\,r\,|P_1 \phi |^2\,d^3 x - \int_{\R^3} \frac{1}{r}\,|\phi|^2\,d^3 x\right)=\nu \left(K[\phi] - \int_{\R^3} \frac{1}{r}\,|\phi|^2\,d^3 x\right)
\end{eqnarray*}
where the kinetic part is defined by $K[\phi]:= \int_{\R^3} t(r)\,r\,|P_1 \phi |^2\,d^3 x$ and satisfies
\begin{eqnarray*}
&&K[\phi]\\
&& =\int_{\R^3} t(r)\,r\left[ |\partial_r \phi|^2 + \frac{1}{r^2}\,|\sigma \!\cdot\! L\,\phi|^2 + |q|^2\,|\phi|^2\right]d^3x \\
&&\qquad +\int_{\R^3} t(r)\Big[-r \,\partial_r \langle \phi, \sigma \!\cdot\! q\,\phi \rangle + r\,\langle \phi, \sigma \!\cdot\! (\partial_r q)\,\phi \rangle-\partial_r \langle \phi, \sigma \!\cdot\! L\,\phi \rangle \Big]\,d^3 x \\
&&\qquad +\int_{\R^3} t(r)\Big[\langle \sigma \!\cdot\! L\,\phi, \sigma \!\cdot\! q\,\phi \rangle +\langle \sigma \!\cdot\! q\,\phi, \sigma \!\cdot\! L\,\phi \rangle\Big] \,d^3x \ .
\end{eqnarray*}
An integration by parts in the $r$ variable yields
\begin{eqnarray*}
&&K[\phi]\\
&& = \int_{\R^3} t(r)\,r\,\left[ |\partial_r \phi|^2 + \frac{1}{r^2}\,|\sigma \!\cdot\! L\,\phi|^2 + |q|^2\,|\phi|^2\right]d^3x\\
&&\; +\int_{\R^3} t(r) \left[4\,\langle \phi, \sigma \!\cdot\! q\,\phi \rangle +\frac{2}{r}\,\langle \phi, \sigma \!\cdot\! L\,\phi \rangle +\langle \sigma \!\cdot\! L\,\phi, \sigma \!\cdot\! q\,\phi \rangle +\langle \sigma \!\cdot\! q\,\phi, \sigma \!\cdot\! L\,\phi \rangle\right] d^3x\\
&&\; +\int_{\R^3} t'(r) \Big[\,r\,\langle \phi, \sigma \!\cdot\! q\,\phi \rangle + \langle \phi, \sigma \!\cdot\! L\,\phi \rangle \,\Big]\,d^3x\ ,
\end{eqnarray*}
where we have also used that $r\,\partial_rq = q$. Consider now the region of integration where $r \leq R$ and denote the corresponding expression by $K_1[\phi]$. There the derivative of $t(r)$ vanishes and hence collecting terms we find
\begin{eqnarray*}
K_1[\phi]=\int_{r \leq R}t\,r \left[|\partial_r \phi|^2 + \left|\left[\frac{1}{r}\,(\sigma \cdot L + 1) +\sigma\cdot q\right]\phi\right|^2- \frac{1}{r^2}\,|\phi|^2\right] d^3 x \\
+\,2 \int_{r \leq R}\langle \phi, \sigma \cdot q\,\phi \rangle \,d^3 x \ .
\end{eqnarray*}

\medskip At this point we have decoupled the derivatives with respect to $r$ from the magnetic field or, to be precise, from $q$. The problem is that the angular momentum is still coupled to the magnetic field. Obviously
$$
\left|\,\Big[\frac{1}{r}\,(\sigma \cdot L + 1) +\sigma\cdot q\Big]\,\phi\right|\,\geq \left|\,\Big|\frac{1}{r}\,(\sigma \cdot L + 1)\,\phi\,\Big| -|\,\sigma\cdot q\,\phi\,| \,\right| \ .
$$
Further, since the eigenvalues of $\sigma \cdot L+1$ are given by $\pm 1$, $\pm 2 \dots $,
$$
\left\|\,(\sigma \cdot L + 1)\,\phi\,\right\|_{L^2(S^2)} \;\geq\quad \left\|\phi\,\right\|_{L^2(S^2)}\ ,
$$
and we have that
$$
\left\|\Big[\frac{1}{r}\,(\sigma \!\cdot\! L + 1) +\sigma\!\cdot\! q\,\Big]\,\phi\,\right\|_{L^2(S^2)}\kern -9pt \geq \frac{1}{r}\,\| \phi \|_{L^2(S^2)} -\Big\|\,|q|\,\phi\,\Big\|_{L^2(S^2)}\kern -3pt \geq \Big[\frac{1}{r}\,-\frac{r}{2}\Big]\,\| \phi\|_{L^2(S^2)}\ ,
$$
{since $|q(r)| \leq r/2$. For $r \leq R \leq \sqrt 2$, the factor $
[\frac{1}{r}\,-\frac{r}{2}]$ is nonnegative. Since $
[\frac{1}{r}\,-\frac{r}{2}]^2-\frac 1{r^2}=\frac{r^2}{4} - 1$ and since $2\,|q(r)|\leq r$, we obtain the lower bound}
$$
K_1[\phi]\geq \int_{r \leq R} r\left[\,|\partial_r \phi|^2 + {\Big[\frac{r^2}{4} - 2\Big]} |\phi|^2\,\right]d^3 x \ .
$$
The function $\frac{r^3}{4} - 2r$ is a decreasing function on the interval $\big[\,0,\sqrt{2}\,\big]$. Hence
$$
K_1[\phi] \ge \int_{r \leq R} r\,|\partial_r \phi|^2\,d^3 x + {\Big[\,\frac {R^3}4-2R\,\Big]} \int_{r \le R} |\phi|^2\,d^3 x\ ,\quad \mbox{if}\ R<\sqrt{2}\ .
$$

\medskip Next, we look at the {contribution to $K[\phi]$ of the} region where $r \geq R$, which we denote by $K_2[\phi]$,
\begin{eqnarray*}
&&K_2[\phi]\\
&& =t \int_{r \geq R}\kern -6pt t(r)\,r\,\left[ |\partial_r \phi|^2 + \frac{1}{r^2}\,|\sigma \!\cdot\! L\,\phi|^2 + |q|^2\,|\phi|^2\right]d^3x\\
&&\; +\int_{r \geq R}\kern -6pt t(r) \left[4\,\langle \phi, \sigma \!\cdot\! q\,\phi \rangle +\frac{2}{r}\,\langle \phi, \sigma \!\cdot\! L\,\phi \rangle +\langle \sigma \!\cdot\! L\,\phi, \sigma \!\cdot\! q\,\phi \rangle +\langle \sigma \!\cdot\! q\,\phi, \sigma \!\cdot\! L\,\phi \rangle\right] d^3x\\
&&\; -\int_{r \geq R}\frac{t(r)}{r}\,\Big[ \,r\,\langle \phi, \sigma \!\cdot\! q\,\phi \rangle + \langle \phi, \sigma \!\cdot\! L\,\phi \rangle \,\Big]\ ,
\end{eqnarray*}
using the fact that $t'=-t/r$. Collecting the terms, we get
\begin{eqnarray*}
&&K_2[\phi]\\&& = \int_{r \geq R}\kern -12pt t(r)\,r\left[ |\partial_r \phi|^2 + \frac{1}{r^2}\,|\sigma \!\cdot\! L\,\phi|^2 + |q|^2\,|\phi|^2\right]d^3x\\
&& +\int_{r \geq R}\kern -12pt t(r) \left[3\,\langle \phi, \sigma \!\cdot\! q\,\phi \rangle +\frac{1}{r}\,\langle \phi, \sigma \!\cdot\! L\,\phi \rangle +\langle \sigma \!\cdot\! L\,\phi, \sigma \!\cdot\! q\,\phi \rangle +\langle \sigma \!\cdot\! q\,\phi, \sigma \!\cdot\! L\,\phi \rangle\right] d^3x \ .
\end{eqnarray*}
This can be rewritten as
\begin{eqnarray*}
K_2[\phi]&&\geq\int_{r \geq R} t(r)\,r\left[ |\partial_r \phi|^2 + \left|\,\Big[\,\frac{1}{r}\,\Big[\sigma \cdot L+\frac 12\Big]+\sigma \cdot q\,\Big] \phi \,\right|^2 -\frac{1}{4\,r^2}\,|\phi|^2\right]d^3x\\
&& \qquad+\;2 \int_{r \geq R} t(r) \,\langle \phi, \sigma \cdot q\,\phi \rangle \,d^3x \ .
\end{eqnarray*}
Finally we get
$$
K_2[\phi] \ge \int_{r \geq R} t(r)\,r\,|\partial_r \phi|^2 \,d^3x - \frac{1}{4R} \int_{r \geq R} |\phi|^2 \,d^3x -2R \int_{r \geq R} \frac{1}{r}\,|q|\,|\phi|^2 \,d^3x\ ,
$$
{and, using $|q|/r\leq 1/2$,}
$$
K_2[\phi] \geq \int_{r \geq R} t(r)\,r\,|\partial_r \phi|^2 \,d^3x - \Big[R+\frac{1}{4R}\Big] \int_{r \geq R} | \phi|^2 \,d^3x\ .
$$

\medskip Thus {we can estimate $K[\phi]=K_1[\phi]+K_2[\phi]$ as follows:}
$$
K[\phi] \geq \int_{\R^3} t(r)\,r\,|\partial_r \phi|^2 \,d^3x + \Big[\frac{R^3}{4}-2R\Big]\int_{r \leq R} |\phi|^2\,d^3 x - \Big[R +\frac{1}{4R}\Big]\int_{r \geq R} | \phi|^2 \,d^3x \ .
$$
Observe that generally
\begin{eqnarray*}\label{posit}
0&\leq& \int_{\R^3}\kern -6pt t(r)\,r\left| \partial_r \phi + \frac{1}{r}\,\phi\right|^2\,d^3 x\\
&&\qquad = \int_{\R^3}\kern -6pt t(r)\,r\,|\partial_r \phi|^2\,d^3 x - \int_{\R^3}\kern -6pt t(r)\,\frac{1}{r}\,|\phi|^2\,d^3 x -\int_{\R^3}\kern -6pt t'(r)\,|\phi|^2 \,d^3x \ .
\end{eqnarray*}
Since $t'(r) \equiv 0$ and $t(r) \equiv 1$ on $[0,R)$, and $t'(r)=-t(r)/r$ on $(R, \infty)$, {the following estimate holds}
\begin{eqnarray*}
\int_{\R^3}\! t(r)\,r\,|\partial_r \phi|^2\,d^3 x - \int_{\R^3} \frac{1}{r}\,|\phi|^2\,d^3 x &\geq& -\int_{\R^3} \!\frac{1-t(r)}{r}\,|\phi|^2\,d^3 x +\int_{\R^3} t'(r)\,|\phi|^2\,d^3 x \\
&&\;= -\int_{r \ge R} \frac{1}{r}\,|\phi|^2 \,d^3x \ge -\frac{1}{R} \int_{r \ge R} |\phi|^2 \,d^3x
\end{eqnarray*}
and hence
$$
K[\phi] - \int_{\R^3} \frac{1}{r}\,|\phi|^2 \,d^3x \geq {\Big[\frac{R^3}{4}-2R\Big]} \int_{r \leq R} |\phi|^2\,d^3 x - \Big[R + \frac{5}{4 R}\Big]\int_{r \geq R} |\phi|^2\,d^3 x \ .
$$
{Optimizing on $R\in (0,\sqrt 2]$, {\sl i.e.\/} using
\[
\max_{R\in(0,\sqrt 2)}\min\left\{\Big[\frac{R^3}{4}-2R\Big],\Big[-R - \frac{5}{4 R}\Big]\right\}=\Big[-R - \frac{5}{4 R}\Big]_{|R=\sqrt 5/2}=-\sqrt 5\ ,
\]
we get
$$
G_1[\phi]\ge\nu \left(K[\phi] - \int_{\R^3} \frac{1}{r}\,|\phi|^2\,d^3 x\right)\ge -\nu\,\sqrt 5 \int_{\R^3} |\phi|^2\,d^3 x \ .
$$
Hence the condition $G_B[\phi] =\sqrt B\,G_1[\phi] \ge -2\,\| \phi \|^2$ entails that
$$
B(\nu)\geq\frac{4}{5\,\nu^2}\ .
$$
}\finprf

\section{Asymptotics for the critical magnetic field}\label{Sec:Asymp}

{In the large magnetic field limit, the upper component of the eigenfunction corresponding to the lowest energy levels in the gap of Dirac operator with magnetic field $H_B-\nu\,|\cdot|^{-1}$ is expected to behave like the eigenfunctions associated to the lowest levels of the Landau operator
\[
L_B:=-\,i\,\sigma_1\,\partial_{x_1}-\,i\,\sigma_2\,\partial_{x_2}-\sigma \cdot { A}_B(x)\;,
\]
which can also we written as $L_B=P_B+\,i\,\sigma_3\,\partial_{x_3}$ or $L_B=-\,i\,(\partial_{x_1}+i\,{Bx_2}/2)\,\sigma_1-\,i\,(\partial_{x_2}-i{B x_1}/2)\,\sigma_2$. The goal of this section is to compare the lowest energy levels of $H_B-\nu\,|\cdot|^{-1}$ with its lowest energy levels on a space generated by the lowest energy levels of $L_B$.} The asymptotic analysis for the small coupling limit $\nu\to 0^+$ is not that simple because {the Landau levels are not stable under the action of the kinetic part of the Dirac Hamiltonian.} The way out is to choose a {representation of $H_B-\nu\,|\cdot|^{-1}$} that diagonalizes the kinetic energy in the Dirac Hamiltonian {and to project both upper and lower components on the lowest Landau levels.}

\subsection{Projection on Landau levels}

{To start with, we observe that
\[
P_B^2=L_B^2-\partial_{x_3}^2
\]
and summarize the basic properties of the lowest energy levels of $L_B$.}
\begin{lemma}\label{Lem:LandauLevels}{\rm[\cite{T}, Section 7.1.3]} The operator $L_B$ in $L^2(\R^2, \C^2)$ has discrete spectrum $\{2n B\,:\,n \in \N\}$, each eigenvalue being infinitely degenerate. Moreover the kernel of this operator, that is, the eigenspace corresponding to the eigenvalue $0$, is the set generated by the $L^2$-normalized functions
$$
\phi_\ell:=\frac{B^{(\ell+1)/2}}{\sqrt{2\,\pi\,{2^\ell\,\ell!}}}\,(x_2+i\,x_1)^\ell\,e^{-{B\,s^2}/{4}}\,\binom 01 \,,\quad\ell\in \N\,,\quad s^2=x_1^2+x_2^2\ .
$$
\end{lemma}

Next we diagonalize the free magnetic Dirac Hamiltonian. First we write it in the form
$$
K_B= \left(\begin{matrix}{\mathbb I} & P_B \\ P_B & -{\mathbb I}\end{matrix}\right)=\sqrt{{\mathbb I}+P_B^2}\; \left(\begin{matrix}R & Q \\ Q & -R\end{matrix}\right)\ ,
$$
where $R$ and $Q$ are operators acting on 2 spinors, given by
$$
R = \frac{1}{\sqrt{{\mathbb I} +P_B^2}} \,, \quad Q = \frac{P_B}{\sqrt{{\mathbb I}+P_B^2}} \ ,
$$
and satisfy the relation
$$
R^2+Q^2={\mathbb I} \ .
$$
The matrix
$$
\left(\begin{matrix}R & Q \\ Q & -R\end{matrix}\right)
$$
is a reflection matrix and hence has eigenvalues $1$ and $-1$. It can be diagonalized using the matrix
$$
U= \frac 1{\sqrt{2\,({\mathbb I}-R)}}\,\left(\begin{matrix}Q &R-{\mathbb I} \\ {\mathbb I}-R& Q \end{matrix}\right)\,.
$$
The operator defined by $U$ is unitary and such that
$$\label{matrixU}
U^*K_B\,U=\left(\begin{matrix} \sqrt{{\mathbb I}+P_B^2} & 0 \\ 0 & -\sqrt{{\mathbb I} +P_B^2}\end{matrix}\right)\,.
$$
The potential $V=\frac1{r}$ is transformed into the nonnegative operator
$$
P:=U^*VU= \left(\begin{matrix}p & q \\ q^* & t\end{matrix}\right) \ .
$$
Here and from now on, we will omit ${\mathbb I}$ whenever it is multiplied by 	a scalar valued function. If we denote by $\W$ any $4$-spinor and decompose it as
$$
\W=\binom \X\Y\ ,
$$
{where $\X$ and $\Y$ are the upper and lower components, in the new representation, the full magnetic Dirac Hamiltonian takes the form
\[
U^*H_B\,U=U^*K_B\,U-\nu\,P=\left(\begin{matrix} \sqrt{{\mathbb I}+P_B^2} & 0 \\ 0 & -\sqrt{{\mathbb I} +P_B^2}\end{matrix}\right)-\nu\,\left(\begin{matrix}p & q \\ q^* & t\end{matrix}\right)\ .
\]
}
The Dirac energy for an electronic wave function $\W$ in the electromagnetic potential $(V,A)$ is now
\[\begin{array}{c}
{\cal E}_\nu[\W]:= \,{\cal K}[\W]-\nu\,(\W, P\,\W )\ , \\ \\
{\cal K}[\W]:=\left(\W, U^* K_B\, U \W\right) = \left(\X, \sqrt{{\mathbb I}+P_B^2}\,\X\right) - \left(\Y, \sqrt{{\mathbb I}+P_B^2}\,\Y\right)\ .
\end{array}\]
As we shall see below, in the new representation, restricting the upper and lower components $\X$ and $\Y$ to the lowest Landau levels makes sense for studying the regime of asymptotically large $B$. The price we pay for that is that all quantities like $R$, $Q$, $U$, $P$... depend on $B$. Denote by $\PiL$ the projector on the lowest Landau level, whose image is generated by all functions
\[
(x_1,x_2,x_3)\mapsto\phi_\ell(x_1,x_2)\,f(x_3)\quad\forall\;\ell\in\N\ ,\quad\forall\;f\in L^2(\R)\ ,
\]
and define $\PiLC:={\mathbb I}-\PiL$. Notice that $\PiL$ commutes with $L_B$. With the above notations, for all $\xi\in L^2(\R^3,\C^2)$, we have that
\be\label{d1}
\left(\xi, \sqrt{{\mathbb I}+P_B^2}\,\xi \right) = \left( \PiL \xi, \sqrt{{\mathbb I}+P_B^2}\,\PiL \xi \right) + \left( \PiLC\xi, \sqrt{{\mathbb I}+P_B^2}\,\PiLC\xi \right)\ .
\ee
Next, we decompose any $\W\in (L^2(\R^3,\C))^4$ as
$$ \W=\PiTot\, \W+\PiTotC\, \W\ ,$$
where
$$
\PiTot := \left(\begin{matrix}\PiL & 0\\ 0&\PiL \end{matrix}\right)\;,\quad \PiTotC := \left(\begin{matrix}\PiLC & 0\\ 0&\PiLC\end{matrix}\right) \ .
$$

\subsection{Main estimates}

{}From \eqn{d1}, it follows that
\[
{\cal K}[\W]={\cal K}[\PiTot\,\W]+{\cal K}[\PiTotC\,\W]\ .
\]
Since the operator $P$ is nonnegative, {we also have
\begin{eqnarray*}
&&\label{d2}(\W, P\,\W) \leq\left(1+ \sqrt\nu\,\right)\,\left(\PiL\,\W, P\,\PiTot\, \W\right)+ \left(1+\frac{1}{\sqrt\nu}\right)\,\left(\PiTotC\, \W, P\,\PiTotC\, \W\right) \ ,\\
&&\label{d3}(\W, P\,\W) \geq\left(1- \sqrt\nu\,\right)\, \left(\PiL\,\W, P\,\PiTot\, \W\right)+ \left(1-\frac{1}{\sqrt\nu}\right)\,\left(\PiTotC\, \W, P\,\PiTotC\, \W\right)\ .
\end{eqnarray*}
This simply follows from the identities
\begin{eqnarray*}
&&(a+b)^2\leq a^2+b^2+2\,|a\,b|=\inf_{\nu>0}\left[\left(1+ \sqrt\nu\,\right)\,a^2+\left(1+\frac{1}{\sqrt\nu}\right)\,b^2\right]\ ,\\
&&(a+b)^2\geq a^2+b^2-2\,|a\,b|=\sup_{\nu>0}\left[\left(1-\sqrt\nu\,\right)\,a^2+\left(1-\frac{1}{\sqrt\nu}\right)\,b^2\right]\ .
\end{eqnarray*}
}
The above remarks prove the next proposition.
\begin{proposition}\label{prop14} For all $\W\in C^\infty_0(\R^3, \C^4]$,
\[{\cal E}_{\nu+\nu^{3/2}}[\PiTot\, \W]+{\cal E}_{\nu+\sqrt\nu}\,[\PiTotC\, \W] \ \leq\ {\cal E}_\nu[\W]\ \leq\ {\cal E}_{\nu-\nu^{3/2}}[\PiTot\, \W]+{\cal E}_{\nu-\sqrt\nu}\,[\PiTotC\, \W]\ .
\]\end{proposition}

The following result will allow us to get rid of the higher Landau levels when looking for the ground state energy, {{\sl i.e.\/} of the term ${\cal E}_{\nu+\sqrt\nu}\,[\PiTotC\, \W]$ in Proposition~\ref{prop14}. Consider $\bar\nu\in (0,1)$ such that
$$2\,(\bar\nu + \sqrt{\bar\nu}\,)=2-\sqrt{2}\ ,$$
{\sl i.e.\/} $\bar\nu\approx 0.056$, and for any $\nu\in (0, \bar\nu)$, define
\[\label{dnu}
d(\delta):=(1-2\delta)\sqrt{2}-2\delta\ ,\quad d_\pm(\nu):=d(\delta_\pm(\nu))\quad\mbox{with}\quad \delta_\pm(\nu):=\sqrt\nu\pm\nu\ .
\]
We have 
\[
d(\delta)>0\quad\Longleftrightarrow\quad\delta <1-\sqrt{2}/2\,,
\]
\[
d_\pm(\nu)>0\quad\mbox{if} \quad\nu<\bar\nu\ .
\]
}
\begin{proposition}\label{prop15} Let $B>0$ and $\delta\in (1-\sqrt{2}/2)$. For any $ \tilde\W=\binom \X0$, $\bar Z=\binom 0\Y$, $\X, \Y \in L^2(\R^3,\C^2)$
\[\label{outHLL}
{\cal E}_\delta\,[\PiTotC\, \tilde\W] \geq d(\delta)\,\sqrt{B}\,\,\|\,\PiLC\,\X\,\|_{_{L^2(\R^3)}}^2\ .
\]
\[\label{outHLL2}
{\cal E}_{-\delta}\,[\PiTotC\, \bar Z] \leq -d(\delta)\,\sqrt{B}\,\,\|\,\PiLC\,\Y\,\|_{_{L^2(\R^3)}}^2\ .
\]
\end{proposition}

\proof {An elementary computation shows that
$$
P_B^2=(\nabla-i{{A}_B})^2\,{\mathbb I}+\sigma\cdot{\mathbf B}\ .
$$
Using the diamagnetic inequality (see \cite{AHS1}),
\[
\int_{\R^3}\Big|\,(\nabla-i{{A}_B})\,\psi\,\Big|^2\,dx\geq \int_{\R^3}\Big|\,\nabla |\psi|\,\Big|^2\,dx\ ,
\]
Hardy's inequality, 
\[
\int_{\R^3}\Big|\,\nabla |\psi|\,\Big|^2\,dx\geq \frac 14\int_{\R^3}\frac{|\psi|^2}{|x|^2}\;dx\ ,
\]
and the nonnegativity of $\sigma\cdot{\mathbf B}+B\,{\mathbb I}$, we get
$$
P_B^2+(1+B)\,{\mathbb I}=(\nabla-i{{A}_B})^2\,{\mathbb I}+\sigma\cdot{\mathbf B}+(1+B)\geq \left(\frac14\,\frac{1}{|x|^2}+1\right)\ .
$$
}
Since the square root is operator monotone, we have
$$
\sqrt B + \sqrt{{\mathbb I}+P_B^2} \ge \sqrt{P_B^2+(1+B)\,{\mathbb I}}\geq \sqrt{\frac14\,\frac{1}{|x|^2}+1}\ .
$$
Now, for any $\delta>0$,
\[\begin{array}{ll}
\sqrt{{\mathbb I}+P_B^2}-\frac{\delta}{|x|}\kern -5pt&= (1-2\delta) \,\sqrt{{\mathbb I}+P_B^2} +2\delta\,\sqrt{{\mathbb I}+P_B^2}-\frac{\delta}{|x|}\\
& \geq (1-2\delta) \,\sqrt{{\mathbb I}+P_B^2}+2\delta \,\sqrt{\frac14\,\frac{1}{|x|^2}+1}- 2\delta\,\sqrt{B}- \frac{\delta}{|x|}\\
& \geq (1-2\delta) \,\sqrt{{\mathbb I}+P_B^2} -2\delta\,\sqrt{B}\ .
\end{array}\]
On the range of $\PiLC$ the operator $\sqrt{{\mathbb I}+P_B^2}\,{\geq\sqrt{{\mathbb I}+L_B^2}}$ is bounded from below by $\sqrt{1+2B}$. {Hence}
$$
\sqrt{{\mathbb I}+P_B^2}-\frac{\delta}{|x|}\ge (1-2\delta) \,\sqrt{1+2B} -2\delta\,\sqrt B\ ,
$$
which is equivalent to
$$
\sqrt{{\mathbb I}+P_B^2} - \delta\,P \ge (1-2\delta) \,\sqrt{1+2B} -2\delta\,\sqrt B \ .
$$
Since
$$
(1-2\delta) \,\sqrt{1+2B} -2\delta\,\sqrt B\sim d(\delta)\,\sqrt B
$$
as $B\to\infty$, and 
$$
\inf_{B>0}[(1-2\delta) \,\sqrt{1+2B} -2\delta\,\sqrt B]>0$$
if $d(\delta)>0$, the right hand side is positive for any field strength provided that $0\le \delta< 1-\sqrt{2}/2$.\finprf

\subsection{The restricted problem}

Next we prove that the ground state energy $\lambda_1(\nu, B)$ is comparable with the one obtained by restricting it to the lowest Landau level, both in the upper and in the lower components of the wave function, provided that the Coulomb potential is slightly modified. By a result similar to Theorem \ref{min-max} (also see Theorem 3.1 in \cite{DES} in case $B=0$), for all $0<B<B(\nu)$, $\nu\in (0,1)$, $\lambda_1(\nu, B)$ is characterized as
\be\label{toto}
\lambda_1(\nu, B)= {\inf_{\aatop{\X\in C^\infty_0(\R^3, \C^2)}{\X\ne 0}}\quad\sup_{\aatop{\Y\in C^\infty_0(\R^3, \C^2)}{\|\W\|_{L^2(\R^3)}=1\,,\;\W=\binom\X\Y}}}{\mathcal E}_\nu[\W]\ .
\ee
With the notation $\W=\binom \X\Y$, we define the restricted
 min-max problem
\[\label{c0}
\llc(\nu,B):= {\inf_{\aatop{\X\in C^\infty_0(\R^3, \C^2)}{\PiLC\X=0\,,\;{0<\|\X\|_{L^2(\R^3)}^2<1}}}\quad\sup_{\aatop{\Y\in C^\infty_0(\R^3, \C^2)\,,\;{\W=\binom\X\Y}}{\PiLC\Y=0\,,\;{\;\|\Y\|_{L^2(\R^3)}=1- \|\X\|_{L^2(\R^3)}^2 }}}} {\mathcal E}_\nu[\W]\ .
\]
We show below that this restricted
 problem is actually a one-dimensional problem. For this purpose, let us define the function $a^B_0:\R\times\R^+\to\R$ given by
$$
a^B_0(z):=B\,\int_0^{+\infty}\frac{s\,e^{-Bs^2/2}}{\sqrt{s^2+z^2}}\,ds\ .
$$
and implicitly define $\mu_{\mathcal L}$ as the unique solution of
$$
\mu_{\mathcal L}[f,\nu,B]=\frac{\displaystyle\int_\R\kern-2pt\left(\kern-1pt\frac{|f'(z)|^2}{1+ \mu_{\mathcal L}[f,\nu,B]+\nu\,a^B_0(z)}+(1-\nu\,a^B_0(z))\,|f(z)|^2\kern-1pt\right)dz}{\displaystyle\int_\R |f(z)|^2\,dz}\ .$$
\begin{theorem}\label{TT} For all $B > 0$ and $\nu\in (0, 1)$,
\[\label{defc0} \llc(\nu,B)= \inf_{f\in C^\infty_0(\R, \C)\setminus\{0\}}\quad \mu_{\mathcal L}[f,\nu,B]\ .\]
\end{theorem}
\proof The definition of $\llc(\nu,B)$ is equivalent to
$$\label{c0phi}
\llc(\nu,B)= {\inf_{\aatop{\phi\in C^\infty_0(\R^3, \C^2)}{\PiLC\phi=0\,,\; \phi\ne 0}}\quad\sup_{\begin{array}{c}\scriptstyle{\chi\in C^\infty_0(\R^3, \C^2)}\cr \scriptstyle{\PiLC\chi=0\,,\; {\psi=\binom\phi\chi}}\cr \scriptstyle\|\psi\|^2_{L^2(\R^3)}=1\end{array}}}\left(\left(\begin{matrix}1-\frac{\nu}{r} & -\,i\,\sigma_3\,\partial_{x_3}\,\\ &\\ -\,i\,\sigma_3\,\partial_{x_3}\, & -1-\frac{\nu}{r} \end{matrix}\right)\!\psi, \psi \right)
$$
or
$$
\llc(\nu,B)= {\inf_{\aatop{\phi\in C^\infty_0(\R^3, \C^2)}{\PiL\phi\ne 0}}\kern-3pt\sup_{\begin{array}{c}{\scriptstyle\chi\in C^\infty_0(\R^3, \C^2)}\cr \scriptstyle{\PiL\chi\neq 0\,,\; {\psi=\binom\phi\chi}}\cr \scriptstyle\|\PiTot\,\psi\|^2_{L^2(\R^3)}=1\end{array}}}\left(\left(\begin{matrix} 1-\frac{\nu}{r} & -\,i\,\sigma_3\,\partial_{x_3}\,\\ &\\ -\,i\,\sigma_3\,\partial_{x_3}\, & -1-\frac{\nu}{r} \end{matrix}\right)\!\PiTot\,\psi, \PiTot\,\psi \right)
$$
with the notation $r=|x|=\sqrt{x_1^2+x_2^2+x_3^2}$. For any given $\phi$ such that $\PiLC\phi=0$, the supremum in $\chi$ is achieved by the function
$$
\chi_{\mathcal L}[\phi]= \Big(V_{\mathcal L}+\big(\lambda_{\mathcal L}[\phi,\nu,B]+1\big)\,\PiL \Big)^{-1} \PiL\,(-\,i\,\sigma_3\,\partial_{x_3})\,\phi\ ,
$$
with ${V_{\mathcal L}(x):=\PiL\,\frac\nu{r}\,\PiL}$ and
$$
\lambda_{\mathcal L}[\phi,\nu,B]:= \kern -7pt\sup_{\aatop{\chi\in C^\infty_0(\R^3, \C^2)}{\PiL\chi\neq 0\,,\; \psi = \binom\phi\chi}}\kern -5pt\left(\left(\begin{matrix}1-\frac{\nu}{r} & -\,i\,\sigma_3\,\partial_{x_3}\,\\ &\\ -\,i\,\sigma_3\,\partial_{x_3}\, & -1-\frac{\nu}{r} \end{matrix}\right)\!\frac{\PiTot\,\psi}{\|\PiTot\,\psi\|_{L^2(\R^3)}}, \frac{\PiTot\,\psi}{\|\PiTot\,\psi\|_{L^2(\R^3)}}\right).
$$
Since $\sigma_3\,\sigma_3^*=\sigma_3^2={\mathbb I}$, this yields the expression
$$
\llc(\nu,B)=\kern -12pt\inf_{\aatop{\phi\in C^\infty_0(\R^3, \C^2)}{\PiLC\phi=0\,,\;\|\phi\|_{L^2(\R^3)}=1}}\int_{\R^3}\left[\frac {|\partial_{x_3}\phi|^2}{V_{\mathcal L}+\lambda_{\mathcal L}[\phi,\nu,B]+1}+\Big(1-\frac{\nu}{r}\Big)|\phi|^2\right]\,dx\ .
$$
Now, with the notations of Lemma~\ref{Lem:LandauLevels}, for all $\ell\ne\ell' \geq 0$, for all $h:\R\to \R$, we have $\int_{\R^2} \phi_\ell \,h(r)\,\phi^*_{\ell'}\,dx_1\,dx_2 =0$, with $s=\sqrt{x_1^2+x_2^2}$, $r=\sqrt{s^2+x_3^2}$, and
$$
\llc(\nu,B)=\inf_{\ell\in \N} \;\inf_{\aatop {\|\phi\|_{L^2(\R^3)}=1}{\phi\in C^\infty_0(\R,\mbox{\tiny span}(\phi_\ell))}}\kern-6pt\Big( (V_{\mathcal L}+\lambda_{\mathcal L}[\phi,\nu,B]+1)^{-1} \partial_{x_3}\phi ,\partial_{x_3}\phi )+\Big(1-\frac{\nu}{r}\Big)\phi, \phi\Big)\ .
$$
A simple calculation shows that for any $g\in C^\infty_0(\R^2,\C^2)$,
$$
\big((V_{\mathcal L})\,g, g\big)_{L^2(\R^2,\C^2)} = \sum_{\ell\geq 0}\,(g, \phi_\ell)_{L^2(\R^2,\C^2)}^2\,\left(\phi_\ell, \frac\nu{r}\,\phi_\ell\right)_{L^2(\R^2,\C^2)} \quad\mbox{a.e. in} \; \R\ni x_3 \ ,
$$
and also that
$$
\left(\phi_\ell, \frac 1r\,\phi_\ell\right)_{L^2(\R^2,\C^2)}=\quad a^B_\ell(x_3)\quad:\,=\quad \frac{B^{\ell+1}}{2^\ell\,\ell!}\,\int_0^{+\infty}\frac{s^{2\ell+1}\,e^{-Bs^2/2}}{\sqrt{s^2+x_3^2}}\,ds \ .
$$
A simple integration by parts shows that for all $\ell\geq 0$, $a^B_\ell\leq a^B_{\ell-1}$ a.e. When minimizing, only the $\ell=0$ component has therefore to be taken into account. \finprf

\begin{corollary}\label{Cor:LandauMonotonicity} For all $\nu\in (0, 1)$, the function $[0,+\infty)\ni B\mapsto \llc(\nu,B)$ is nonincreasing in $B$. \end{corollary}
\proof A simple change of variables shows that
$$
a^B_0(z)=\,\int_0^{+\infty}\frac{s\,e^{-s^2/2}}{\sqrt{\frac{s^2}B+z^2}}\,ds\ .
$$
By Theorem~\ref{TT} and according to the definition of $\mu_{\mathcal L}[f,\nu,B]$, this implies the monotonicity of $\llc(\nu,\cdot)$ in $[0,\infty)$.\finprf

\begin{proposition}\label{prop17} For all $B\geq 0$, the function $\nu\mapsto \llc(\nu,B)$ is continuous in the interval $(0, 1)$ as long as it takes its values in $(-1,1)$. {Moreover, for any $\nu\in (0,1)$, as long as $\llc(\nu,B)$ takes its values in $(-1,1)$, there exists a function $\W\in\mbox{Range}(\PiTot)$ with $\|\W\|_{L^2(\R^3,\C^4)}=1$ such that ${\mathcal E}_\nu[\W]=\llc(\nu,B)$.} \end{proposition}
The proof is similar to the one of Proposition \ref{cont}.\finprf

The above proposition enables us to define
$$\BL(\nu) := \inf\{B>0\;:\; \llc(\nu,B)=-1\}\ .$$
Recall that $\lim_{\nu\to 0^+}d_+(\nu)^{-2}=1/2$. We are now ready to state and prove the main result of this section.

\subsection{Asymptotic results}

\begin{theorem}\label{thmfirstll} Let $\nu\in(0,\bar\nu)$. For any $B\in\big(1/d_+(\nu)^{2},\min\big\{B(\nu), \BL(\nu+\nu^{3/2})\big\}\big)$, we have
\[\label{firstfirst}
\llc\Big(\nu+\nu^{3/2}, B\Big)\ \leq\ \lambda_1(\nu, B)\ \leq\ \llc(\nu-\nu^{3/2},B)\ .
\]\end{theorem}
Notice that the right hand side inequality holds for any $B\in\big(1/d_+(\nu)^{2},B(\nu)\big)$. 

\proof To prove the upper estimate, we use \eqn{toto} and notice that 
\[\label{minoriz}
\lambda_1(\nu, B)\leq \inf_{\aatop {X\in C^\infty_0(\R^3, \C^2)}{\PiLC X=0, \PiL X\ne 0}}\;\;\sup_{\aatop{Y\in C^\infty_0(\R^3, \C^2)}{ ||\W||_{L^2(\R^3)}=1\,,\;\W=\binom\X\Y}}{\mathcal E}_\nu(\W)\,,
\]
since adding the condition $\PiLC X=0$ increases the value of the infimum. Then, by Propositions \ref{prop14} and \ref{prop15},
\begin{eqnarray*}
\lambda_1(\nu,B) &\leq& \inf_{\aatop{\X\in C^\infty_0(\R^3, \C^2)}{ \PiLC \X=0, \PiL \X\ne 0}}\sup_{{\Y\in C^\infty_0(\R^3, \C^2)}}{\cal E}_{\nu-\nu^{3/2}}[\PiTot\, \W]+{\cal E}_{\nu-\sqrt\nu}\,[\PiTotC\, \W]\\
&\leq& \inf_{\aatop{\X\in C^\infty_0(\R^3, \C^2)}{ \PiLC \X=0, \PiL \X\ne 0}}\sup_{{\Y\in C^\infty_0(\R^3, \C^2)}}\frac{{\cal E}_{\nu-\nu^{3/2}}[\PiTot\,\W]-d_-(\nu)\,\sqrt{B}\,\|\PiLC\Y\|_{_{L^2(\R^3)}}^2}{\|\PiTot\,\W\|_{L^2(\R^3)}^2+ \|\PiLC\Y\|_{L^2(\R^3)}^2}\\
&\leq& \inf_{\aatop{\X\in C^\infty_0(\R^3, \C^2)}{ \PiLC \X=0, \PiL \X\ne 0}}
\sup_{{\Y\in C^\infty_0(\R^3, \C^2)}}
\frac{{\cal E}_{\nu-\nu^{3/2}}[\PiTot\,\W]}{\|\PiTot\,\W\|_{L^2(\R^3)}^2}= \llc(\nu-\nu^{3/2},B)\,.
\end{eqnarray*}

\medskip Next, we establish the lower bound. By taking in (\ref{toto}) a smaller maximizing class of functions we decrease the maximum:
$$
\lambda_1( \nu, B)\geq {\inf_{\aatop{\X\in C^\infty_0(\R^3, \C^2)}{\X\ne 0}}\quad\sup_{\aatop{\Y\in C^\infty_0(\R^3, \C^2)\,,\;\W=\binom\X\Y}{\PiLC\Y=0\,,\;\|\W\|_{L^2(\R^3)}=1}}}{\mathcal E}_\nu[\W]\ .
$$
Therefore, by Propositions \ref{prop14} and \ref{prop15}, 
\begin{eqnarray*}
\lambda_1(\nu,B)&\geq& {\inf_{\aatop{\X\in C^\infty_0(\R^3, \C^2)}{\X\ne 0}}\quad\sup_{\aatop{\Y\in C^\infty_0(\R^3, \C^2)\,,\;\W=\binom\X\Y}{\PiLC\Y=0\,,\;\|\W\|_{L^2(\R^3)}=1}}}\Big(\; {\cal E}_{\nu+\nu^{3/2}}[\PiTot\,\W]+ {\cal E}_{\nu+\sqrt\nu}[\PiTotC\,\W]\;\Big)\ ,\\
\label{second} &\geq& \inf_{\aatop{\X\in C^\infty_0(\R^3, \C^2)}{\X\ne 0}}\sup_{\aatop{\Y\in C^\infty_0(\R^3, \C^2)\,,\;\W=\binom\X\Y}{\PiLC\Y=0\,,\;\|\W\|_{L^2(\R^3)}=1}}\kern -10pt\frac{{\cal E}_{\nu+\nu^{3/2}}[\PiTot\,\W]+\dnu\,\sqrt{B}\,\|\PiLC\X\|_{_{L^2(\R^3)}}^2}{\|\PiTot\,\W\|_{L^2(\R^3)}^2+ \|\PiLC\X\|_{L^2(\R^3)}^2}\ .
\end{eqnarray*}
Let us now notice that for every $\X\in C^\infty_0(\R^3, \C^2),\; \X\ne 0$, 
$$
\sup_{\aatop{\Y\in C^\infty_0(\R^3, \C^2)}{{\W=\binom\X\Y}}}\frac{{\cal E}_{\nu+\nu^{3/2}}[\PiTot\,\W]}{\|\PiTot\,\W\|_{L^2(\R^3)}^2}=\lambda_{\mathcal L}[\X,\nu+\nu^{3/2},B]
$$
is uniquely achieved at some $\Y_{\mathcal L}[\X]$ because of the same concavity argument as in the proof of Theorem~\ref{min-max}, after noticing that for any $B\in(0, \BL(\nu+\nu^{3/2}))$, $-1<\llc(\nu+\nu^{3/2},B)\leq{\lambda_{\mathcal L}[\X,\nu+\nu^{3/2},B]}$ and $V_{\mathcal L}\geq 0$. Recall that
$$
\llc(\nu+\nu^{3/2}, B) = \inf_{\aatop{\X\in C^\infty_0(\R^3, \C^2)}{\X\ne 0}} \lambda_{\mathcal L}[\X,\nu+\nu^{3/2},B]\ .
$$
Denoting $\W_{\mathcal L}[\X]= \binom \X{\Y_{\mathcal L}[\X]}$, for any given $\X$, we find
\begin{eqnarray*}
&&\kern-48pt\sup_{{\Y\in C^\infty_0(\R^3, \C^2)}}\frac{{\cal E}_{\nu+\nu^{3/2}}[\PiTot\,\W]+\dnu\,\sqrt{B}\,\|\PiLC\X\|_{_{L^2(\R^3)}}^2}{\|\PiTot\,\W\|_{L^2(\R^3)}^2+ \|\PiLC\X\|_{L^2(\R^3)}^2}\\
&\geq& \frac{{\cal E}_{\nu+\nu^{3/2}}(\PiTot\, \W_{\mathcal L}[\X])+\dnu\,\sqrt{B}\,\|\PiLC\X\|_{_{L^2(\R^3)}}^2}{\|\PiTot\,\W_{\mathcal L}[\X]\|_{L^2(\R^3)}^2+ \|\PiLC\X\|_{L^2(\R^3)}^2}\\
&&\kern24pt=\frac{\lambda_{\mathcal L}\X,\nu+\nu^{3/2},B]\|\PiTot\,\W_{\mathcal L}[\X]\|_{L^2(\R^3)}^2+\dnu\,\sqrt{B}\,\|\PiLC\X\|_{_{L^2(\R^3)}}^2}{\|\PiTot\,\W_{\mathcal L}[\X]\|_{L^2(\R^3)}^2+ \|\PiLC\X\|_{L^2(\R^3)}^2}\\
&&\kern48pt\geq\llc(\nu+\nu^{3/2},B)
\end{eqnarray*}
for $B$ large enough so that $\dnu\sqrt{B}\geq\llc(\nu+\nu^{3/2},B)$. As we shall see below, this is always possible.

Note indeed that on $(0, \bar \nu)$, $d_+(\nu)\leq \sqrt2$. Hence, $d_+(\nu)^{-2}\geq 1/2$. Now, by monotonicity (see Corollary~\ref{Cor:LandauMonotonicity}), $\llc(\delta,B)\leq \llc(\delta,1/2)$ for all $\delta\in (0,\bar\nu+\sqrt{\bar\nu})=(0,1-\sqrt 2/2)$ and for all $B\geq 1/2$. Moreover, one can prove very easily that $\llc(\delta,1/2)\leq 1$ for all $\delta \in (0,1)$.
Indeed, by Theorem \ref{TT}, for all $B$, $\llc(\delta,B) =\inf_{f} \mu_{\mathcal L}[f,\nu,B]$. A simple scaling argument shows that we can make $\int_{\R}|f'|^2\,dz$ as small as we wish while keeping $\int_{\R}|f|^2\,dz$ constant. Taking into account the definition of $\mu_{\mathcal L}[f,\nu,B]$, this shows that $\,\llc(\delta,B) \leq 1$ for all $B$, for all $\delta\in (0,1)$. Therefore, for all $B\geq d_+(\nu)^{-2}$, $\dnu\sqrt{B}\geq\llc(\nu+\nu^{3/2},B)$ holds true.
\finprf

{}From Theorem~\ref{thmfirstll}, we deduce the following
\begin{corollary}\label{Cor:BL} Let $\nu\in (0,\bar\nu)$. Then
\be\label{LMLM} \BL(\nu+\nu^{3/2})\ \leq\ B(\nu)\ \leq\ \BL(\nu-\nu^{3/2})\ .\ee
\end{corollary}

{}{From Theorem~\ref{TT} and Corollary \ref{Cor:BL}, better estimates of the critical magnetic strength $B(\nu)$ than those of Theorem \ref{HHH} can be established for $\nu$ small.}
\begin{theorem}\label{CC} The critical strength $B(\nu)$ satisfies:
$$ \lim_{\nu \to 0} \nu \log B(\nu)\,= \pi\ .$$
\end{theorem}

\proof {Because of Corollary~\ref{Cor:BL}, $B(\nu)$ can be estimated using $\BL(\deltanu)$ with $\deltanu=\nu\pm\nu^{3/2}$. This amounts to look for the smallest positive $B$ for which
\[\label{Eqn:LandauScaling}
\lambda_{\mathcal L}(\deltanu, B):=1+\inf_{\aatop{f\in C^\infty_0(\R, \C)\setminus\{0\}}{\|f\|_{L^2(\R)}=1}}\quad \int_\R\left(\frac{|f'(z)|^2}{\deltanu\,a^{B}_0(z)}-\deltanu\,a^B_0(z)\,|f(z)|^2\right)\,dz
\]
is such that}
\[\label{A0}
\lambda_{\mathcal L}(\deltanu, B)= -1\ .
\]
{Using the identity
\[
a_0^B(z)=\sqrt B\,a^1_0\left(\sqrt B\,z\right)\ ,
\]
}by the changes of variable and function
$$
y(z) := \int_0^z a^1_0(t)\,dt\ ,\quad f\left({\frac z{\sqrt B}}\right)= B^{1/4}\,g(y)\ ,
$$
one tranforms the {above minimization problem into}
\[
\lambda_{\mathcal L}(\deltanu, B)-1=\sqrt{B}\,\inf_{\aatop{g\in C^\infty_0(\R, \C)\setminus\{0\}} {\int_\R g(y)^2\,{d\mu(y)}=1}}\quad \int_\R\left(\frac 1\deltanu\,|g'(y)|^2-\deltanu\,|g(y)|^2\right)\,dy\ .
\]
Hence, 
\be\label{L0}
\lambda_{\mathcal L}(\deltanu, B)=1+\sqrt{B}\left(\lambda_{\mathcal L}(\deltanu, 1)-1\right)\ .
\ee
For a given $\deltanu$, let $\kappa=\kappa(\deltanu):=\deltanu\,\big(1-\lambda_{\mathcal L}(\deltanu, 1)\big)$ and $\mu(y):=1/a_0^1(z(y))$. The problem is reduced to look for the first eigenvalue $E_1=E_1(\deltanu)$ of the operator $-\partial_y^2+\kappa(\deltanu)\,\mu(y)$, namely to find $\deltanu$ such that
\[\label{eqqq}
\deltanu^2=E_1(\delta)\ .
\]
The function $a^1_0$ satisfies
$$ a^1_0(z)\leq a^1_0(0)=\sqrt{\frac{\pi}{2}} \quad\forall\;z\in\R\ , \quad a^1_0(z)\sim \frac{1}{|z|}\;\mbox{as}\;|z|\to\infty\ .
$$
{There exists therefore a constant $c>0$ such that $\mu(y)\leq c\, e^{|y|}$ for any $y\in\R$.} To get an upper estimate of $E_1$, we may now consider the function $g_1(y):=\cos (\pi\,y/2)$ {on $(-1,1)\ni y$} and the rescaled functions $g_\sigma(y):=\sigma^{-1/2}g_1\big(\cdot/\sigma)$.
\begin{eqnarray*}
E_1(\delta)\leq \frac{\pi^2}{4\,\sigma^2}+ \kappa \int_{-1}^1 |g_1|^2\,\mu(\sigma\,y)\,dy&\leq&\frac{\pi^2}{4\,\sigma^2}+ \kappa\,c\int_{-1}^1e^{\sigma|y|}|g_1|^2\,dy\\
&&\qquad\leq \frac{\pi^2}{4\,\sigma^2}+ 2\,\kappa\,c\left(e^{\sigma}-1\right)\ .
\end{eqnarray*}
Optimizing in $\sigma$ in the above expression, we choose $\sigma=\sigma(\deltanu)$ satisfying
$$
\pi^2=4\,\kappa\,c\,\sigma^3\,e^{\sigma}\ ,
$$
which implies
$$
\sigma(\deltanu)\sim -\log\kappa=:\sigma_\deltanu\to \infty\quad\mbox{as}\quad\deltanu\to 0\ .
$$
A Taylor expansion at next order shows that 
\[
\sigma(\deltanu)-\sigma_\deltanu\sim -3\log\big(\sigma_\deltanu\big)\ ,
\]
which yields
$$
E_1(\delta)\leq\frac{\pi^2}{4\,\sigma_\deltanu^2}\,(1+o(1))\ .
$$

\medskip Next, in order to obtain a lower estimate of $E_1(\delta)$, we consider the function $\mu_{\deltanu}$ which is equal to $0$ in the interval $(-\sigma_\deltanu, \sigma_\deltanu)$ and equal to $\kappa\,\mu(\sigma_\deltanu)$ elsewhere. The function $\mu_{\deltanu}$ is positive, but below the function $ \kappa\,\mu(y)$. Then observe that the first eigenvalue of the operator $-\partial_y^2 +\mu_{\deltanu}(y) $, that we denote $E_1^{\deltanu}$, satisfies the equation
$$
\sqrt{E_1^{\deltanu}}\,\sigma_\deltanu = \arctan\left(\sqrt{\frac{{\kappa\,\mu(\sigma_\deltanu)-E_1^{\deltanu}}}{{E_1^{\deltanu}}}}\,\right)\ ,
$$
and as $ \deltanu$ goes to $0^+$ this implies
$$
E_1(\delta)\geq E_1^{\deltanu}=\frac{\pi^2}{4\,\sigma_\deltanu^2}\,\big(1+o(1)\big)\,.
$$
Summarizing, what we have obtained is
\[
E_1(\delta)=\frac{\pi^2}{4\,\big(\log\kappa(\deltanu)\big)^2}\,\big(1+o(1)\big)\,.
\]
So, imposing $E_1(\delta)=\deltanu^2$, we get $\kappa(\delta)=e^{-\frac{\pi}{2\delta}(1+o(1))}$. Since by \eqref{L0}, 
\[
-1=\lambda_{\mathcal L}\big(\deltanu, \BL(\delta)\big)=1+\sqrt{\BL(\delta)}\,\,\frac{\kappa(\deltanu)}\deltanu\,,
\]
we get
$$\BL(\deltanu)\,=\,4\,\deltanu^2\,e^{\frac\pi\deltanu\,(1+o(1))}\ ,$$
which, together with \eqref{LMLM}, concludes the proof.\finprf

\section*{Appendix\footnote{The following proof was explained to us by George Nenciu to whom we are grateful.}: selfadjointness of $H_B-\nu\,|\cdot|^{-1}$}

With the notation $P=-i\,\nabla$, consider {a} Dirac operator of the form
\[\label{H}
{H_0}={\alpha \cdot P} +m\,{\beta} + {\mathbb{V}_0}({x})
\]
defined on $(\mathcal{C}_{0}^{\infty}(\mathbb{R}^{3}))^{4}$. {We assume for instance} that ${(\mathbb{V}_0)_{i,j}}({x})\in L^{2}_{\rm loc}(\mathbb{R}^{3})$. If $f\in \mathcal{C}_{0}^{\infty}(\mathbb{R}^{3})$, {then} the following identity holds on $(\mathcal{C}_{0}^{\infty}(\mathbb{R}^{3}))^{4}$ :
\begin{equation}\label{id}
{H_0}\,f-f\,{H_0}-i\,{\alpha} \cdot \nabla\,f=0 \ .
\end{equation}
{Still denote by} $H_{0}$ the a selfadjoint extension of ${\alpha \cdot P} +m\,{\beta} +\mathbb{V}_{0}({x})$ with domain $\mathcal{D}_{0}$ and let $f\in \mathcal{C}_{0}^{\infty}(\mathbb{R}^{3})$ be such that $0\leq f({x}) \leq 1$, $f({x}) =1$ for $|{x} |\leq 1$, $f({x}) =0$ for $|{x} |\geq 2$, and {$f_a({x}):= f({x}/a)$.} If
\[
\mathcal{D}:=\Big\{g\equiv f_a\,\psi\;|\;a\geq 1,\;\psi \in \mathcal{D}_{0} \Big\}\;,
\]
then in all interesting cases, including the case of local Coulomb singularities with $\nu < 1$ (see the characterisation of $\mathcal{D}_{0}$ in \cite{KW}), one has $\mathcal{D} \subset \mathcal{D}_{0}$, and then by density, Identity~(\ref{id}) {implies that} for all $\psi \in \mathcal{D}_{0}$,
\begin{equation}\label{id1}
\big(H_{0}\,f_a-f_a\,H_{0}-i\,{\alpha} \cdot \nabla\,f_a\big)\,\psi=0 \ .
\end{equation}

\noindent{
A standard characterization of selfadjoint operators that we are going to use is the following: {\sl Let $T$ be a closed symmetric operator. $T$ is selfadjoint if and only if\/}
\[{\rm Ker}\,(T^*\pm\,i)=\{0\}\;.\]}

The first remark is that $H_{0}$ is essentially selfadjoint on $ \mathcal{D}$. Suppose indeed that
\[
\langle\Psi, (H_{0}\pm\,i)\,{(f_a\psi)}\rangle=0\quad {\forall\; a \geq 1\,,\quad\forall\;\psi \in \mathcal{D}_{0}}
\]
{for some $\Psi\in(L^2(\R^3))^4$.} Now from (\ref{id1}),
$$
\langle\Psi, (H_{0}\pm\,i)\,f_a\psi\rangle= \langle\Psi, (H_{0}\pm\,i)\,\psi\rangle -\langle(1-f_a)\Psi, (H_{0}\pm\,i)\,\psi\rangle+i\,\langle\Psi,({\alpha \cdot \nabla}f_a)\psi \rangle \ ,
$$
so that, taking $a\rightarrow \infty$ at fixed $\psi$, one obtains
\[
\langle\Psi, (H_{0}\pm\,i)\,\psi\rangle=0\quad {\forall\;\psi \in \mathcal{D}_{0}\ ,}
\]
and then $\Psi =0$.

{This result also applies} to the case
$$
\mathbb{V}_{0} = -\frac{\nu}{|x|}\,,\quad 0 \le \nu < 1 \ .
$$
{See \cite{KW,N} for more details. Define $\nnrm\cdot$ as the matrix norm. Suppose that $\mathbb{V}_{1}$ is locally~$L^{\infty}$ and more precisely satisfies}
\[\label{loc}
M(R):=\sup_{\vert {x} \vert \leq R}\nnrm{\mathbb{V}_{1}({x})} < \infty\quad\forall\; R\in\R^+\,, \quad \lim_{R\rightarrow \infty}M(R)=\infty\ .
\]
Consider on $\mathcal{D}$ the operator $ H=H_{0}+\mathbb{V}_{1}$.
\begin{lemma} Under the above assumptions, $H$ is essentially selfadjoint on $\mathcal{D}$. \end{lemma}
\proof Let $\chi_{R}$ be the characteristic function of the set $\{ {x} \;|\; \vert {x} \vert \leq {2}R \}$. Since $\chi_{R}\mathbb{V}_{1}$ is bounded, $H_{R}=H_{0}+\chi_{R}\mathbb{V}_{1}$ is essentially selfadjoint on $\mathcal{D}$. {Suppose that there exists $\Psi\in (L^2(\R^3))^4$ such that}
\[\label{con}\langle\Psi, (H_{0}{\pm\,i})\,{(f_a\psi)}\rangle=0\quad {\forall\; a \geq 1\,,\quad\forall\;\psi \in \mathcal{D}_{0}}
\]
{and assume that $\|\Psi\|_{L^2(\R^3)} =1$. Using} (\ref{id}), one deduces that
\[
<f_a\Psi, (H {\pm\,i})\,\psi\rangle ={-i}\ \langle\Psi, {\alpha \cdot \nabla}f_a\psi\rangle\ .
\]
Observe now that $<f_a\Psi, (H {\pm\,i})\,\psi\rangle=\langle f_a\Psi, (H_a{\pm\,i})\,\psi\rangle$, which amounts to 
\begin{equation}\label{com}
\langle f_a\Psi, (H_a{\pm\,i})\,\psi\rangle ={-i}\ \langle\Psi, {\alpha \cdot \nabla}f_a\psi\rangle\ .
\end{equation}
Since $H_a$ is essentially selfadjoint on $\mathcal{D}$, $\mbox{Range}(H_a{\pm\,i})$ is dense in $\mathcal{D}$ and there exists ${\psi_a^\pm} \in \mathcal{D}$ such that
\begin{equation}\label{ex}
(H_a{\pm\,i})\,{\psi_a^\pm}=f_a\Psi +\delta_a
\end{equation}
with $\| \delta_a\|_{L^2(\R^3)} \leq 1/a$. Also notice that
\[\label{symm}
\| (H_a{\pm\,i})\,{\psi_a^\pm} \|_{L^2(\R^3)}^2=\| H_a\,{\psi_a^\pm} \|_{L^2(\R^3)}^2 + \| \psi_a^\pm \|_{L^2(\R^3)}^2 \geq \| {\psi_a^\pm} \|_{L^2(\R^3)}^2\ .
\]

From (\ref{com}) written for ${\psi_a^\pm}$ and (\ref{ex}), {we get}
$$
\| f_a\Psi\|_{L^2(\R^3)} ^{2}+\langle f_a\Psi, \delta_a\rangle= -i\,\langle\Psi, {\alpha \cdot \nabla}f_a{\psi_a^\pm}\rangle
$$
so that
\begin{eqnarray*}
\| f_a\Psi\|_{L^2(\R^3)} ^{2}&\leq& \frac 1a +\|\nabla f_a\|_{L^\infty(\R^3)} \,\| {\psi_a^\pm} \|_{L^2(\R^3)}\\
&\leq& \frac 1a +\|\nabla f_a\|_{L^\infty(\R^3)} \,\| (H_a{\pm\,i})\,{\psi_a^\pm} \|_{L^2(\R^3)}\\
&\leq& \frac 1a +\|\nabla f_a\|_{L^\infty(\R^3)} \,\Big(1+\frac 1a\Big)\ .
\end{eqnarray*}
For $a \rightarrow \infty$, $\| f_a\Psi\|_{L^2(\R^3)} \rightarrow 1$ and $\|\nabla f_a\|_{L^\infty(\R^3)} \rightarrow 0$, a contradiction: $\Psi =0$.\finprf

\bigskip\noindent{\bf Acknowledgments.} We are grateful to George Nenciu for showing us the proof of self-adjointness in the appendix, to P. Pickl and D. D\"urr for references. Some of this research has been carried out at the Erwin Schr\"odinger Institute. M.L. would like to thank the Ceremade for its hospitality. J.D. and M.J.E. acknowledge support from ANR Accquarel project and European Program ``Analysis and Quantum'' HPRN-CT \# 2002-00277. M.L. is partially supported by U.S. National Science Foundation grant DMS 03-00349.

\bigskip\noindent{\sl\scriptsize \copyright~2006 by the authors. This paper may be reproduced, in its entirety, for non-commercial purposes.}

\end{document}